\documentclass[10pt]{amsart}
\usepackage{amscd}
\usepackage{amssymb}
\usepackage[matrix,arrow]{xy}
\usepackage{color}
\newtheorem{theorem}{Theorem}[section]
\newtheorem{proposition}[theorem]{Proposition}
\newtheorem{lemma}[theorem]{Lemma}

\newtheorem{corollary}[theorem]{Corollary}
\theoremstyle{remark}
\newtheorem{remark}[theorem]{Remark}
\newtheorem{examples}[theorem]{Examples}
\newtheorem{definition}[theorem]{Definition}

\newcommand{\id}{\mathrm{id}}

\newcommand{\End}{\mathrm{End}}
\newcommand{\Aut}{\mathrm{Aut}}

\newcommand{\ot}{\otimes}

\newcommand{\SSS}{\scriptscriptstyle}

\newcommand{\U}{\mathcal{U}}
\newcommand{\ho}{homomorphism}
\newcommand{\mo}{monomorphism}
\newcommand{\hos}{homomorphisms}
\newcommand{\mos}{monomorphisms}

\newcommand{\CCC}{\mathbb{C}}

\newcommand{\ep}{\varepsilon}

\newcommand{\ZZZ}{\mathbb{Z}}
\newcommand{\Z}{\mathbb{Z}}

\newcommand{\OOO}{\mathcal{O}}

\newcommand{\weakr}{weakly rigid}

\newcommand{\fset}{\mathcal{F}}
\newcommand{\hset}{\mathcal{H}}

\newcommand{\gset}{\mathcal{G}}

\newcommand{\wsp}{weakly semiprojective}

\newcounter{rocount}

\begin{document}
\title{Continuous fields of C*-algebras over finite dimensional spaces}
\thanks{The author was partially supported by
NSF grant \#DMS-0500693}

\maketitle

 \centerline{\normalsize
MARIUS DADARLAT} \vskip 4pt  \centerline{\it \footnotesize Purdue University}
\centerline{\it\footnotesize West Lafayette, IN, U.S.A.} \vskip 4pt

%\author{Marius Dadarlat}
%
%\address{Department of Mathematics, Purdue University, West
%Lafayette IN 47907, U.S.A.} \email{mdd@math.purdue.edu}
%\date{\today}

%\subjclass{46L07, 47A99}

\begin{abstract} Let $X$ be a finite dimensional compact metrizable space.
We study a technique which employs semiprojectivity  as a tool to produce
approximations of $C(X)$-algebras by $C(X)$-subalgebras with controlled
complexity.
The following applications are given.
All unital separable
continuous fields  of C*-algebras over $X$ with fibers isomorphic
 to a fixed Cuntz algebra $\mathcal{O}_n$, $n\in\{2,3,...,\infty\}$
are locally trivial. They are  trivial if $n=2$ or $n=\infty$. For  $n$ $\geq
3$ finite,
  such a field  is trivial if and only if $(n-1)[1_A]=0$ in $K_0(A)$,
   where $A$ is the C*-algebra of continuous sections of the field.
   We give a complete list of the Kirchberg algebras $D$ satisfying  the
   UCT and having finitely generated K-theory groups for which  every
   unital separable
continuous field  over $X$ with fibers isomorphic to $D$ is automatically
locally trivial or trivial. In a more general context,  we show that  a separable unital
continuous field over $X$ with fibers isomorphic to a $KK$-semiprojective
  Kirchberg C*-algebra is  trivial if and only if it satisfies
  a K-theoretical Fell type condition.
    \end{abstract}

\tableofcontents

\section{Introduction}
 Gelfand's
characterization of commutative C*-algebras
 has suggested the problem of representing non-commutative C*-algebras
as sections of bundles.
 By
a result of Fell \cite{Fell}, if the primitive spectrum  $X$ of a separable C*-algebra $A$ is Hausdorff, then
 $A$ is isomorphic
 to the C*-algebra of continuous sections
vanishing at infinity of a continuous field of simple C*-algebras over $X$.
 In particular $A$ is a continuous $C(X)$-algebra
in the sense of Kasparov \cite{Kas:inv}. This description is very satisfactory,
since as explained in \cite{BK:bundles}, the continuous fields of
C*-algebras are in natural correspondence with the bundles of C*-algebras in
the sense of topology. Nevertheless,
 only a tiny fraction of
 the continuous fields of C*-algebras
 correspond to locally trivial bundles.

In this paper we prove  automatic and conditional local/global trivialization
results for continuous fields of Kirchberg algebras. By a Kirchberg algebra we
mean a purely infinite simple nuclear separable C*-algebra
\cite{Ror:encyclopedia}. Notable examples include the simple Cuntz-Krieger
algebras \cite{CuntzKri:OA}. The following theorem  illustrates   our
results.
\begin{theorem}\label{cuntz-algebras-intro}
    A  separable
unital $C(X)$-algebra $A$ over a finite dimensional
 compact Hausdorff space $X$ all of whose
fibers  are isomorphic
 to the same Cuntz algebra $\mathcal{O}_n$, $n\in\{2,3,\dots,\infty\}$,
is locally trivial. If $\,n=2$ or $n=\infty$, then $A\cong C(X)\ot \OOO_n$.
 If $\,3
\leq n<\infty$, then $A$ is isomorphic to $C(X)\ot\OOO_n$ if and only if
$\,(n-1)[1_A]=0$ in $K_0(A)$.
\end{theorem} The case $X=[0,1]$ of Theorem~\ref{cuntz-algebras-intro}
was proved in a joint paper with G. Elliott
 \cite{DadEll:bundles-interval}.

 We  parametrize the homotopy classes
\[[X,\mathrm{Aut}(\OOO_n)]\cong
 \left\{
\begin{array}{ll}
    K_1(C(X)\ot \OOO_n) &\hbox{if $3\leq n <\infty$}, \\
    \{*\} & \hbox{if $n=2,\,\infty$}, \\
\end{array}
\right.\]  (see Theorem~\ref{nontrivialO(n)bundle}) and hence classify the
unital separable $C(SX)$-algebras $A$ with fiber
 $\OOO_n$ over the suspension $SX$  of a finite dimensional metrizable Hausdorff space $X$.

To put our results in  perspective, let us recall that none of the general
basic properties of a continuous field
 implies any kind of local triviality.
 An example of a continuous field of
 Kirchberg algebras over $[0,1]$
  which is not locally trivial at any point even
though all of its fibers are mutually isomorphic is exhibited in
\cite[Ex.~8.4]{DadEll:bundles-interval}. Examples of nonexact continuous fields
with similar properties  were found by S. Wassermann ~\cite{Wass:to-mdd}.

 A separable C*-algebra $D$ is $KK$-semiprojective if the
functor $KK(D,-)$ is continuous, see Sec.~\ref{section:Semiprojectivity}. The class of $KK$-semiprojective C*-algebras
includes the
 nuclear semiprojective C*-algebras and also the C*-algebras
 which satisfy the Universal
 Coefficient Theorem in KK-theory (abbreviated UCT \cite{RosSho:UCT})
  and whose K-theory groups are finitely generated.
  It is very interesting
that the \emph{only obstruction} to local or global triviality
 for a continuous field of
Kirchberg algebras  is of purely
 K-theoretical nature.
\begin{theorem}\label{stable-trivial-C(X)-algebras}
  Let $A$ be a  separable C*-algebra whose primitive
 spectrum $X$ is compact Hausdorff  and of finite dimension.
 Suppose that each primitive quotient $A(x)$ of $A$ is
 nuclear, purely infinite and stable.
  Then $A$ is isomorphic to $C(X)\ot D$ for some $KK$-semiprojective
stable  Kirchberg algebra $D$ if and only if there is $\sigma \in KK(D,A)$ such
that $\sigma_x \in KK(D, A(x))^{-1}$ for all $x \in X$. For any such $\sigma$
there is an  isomorphism of $C(X)$-algebras $\Phi:C(X)\ot D\to A$ such that
$KK(\Phi|_D)=\sigma$.
 \end{theorem}

  We
have an entirely  similar result  covering the  unital case:
Theorem~\ref{trivial-C(X)-algebras}. The required existence of $\sigma$  is
 a KK-theoretical analog of the classical  condition of Fell
that appears in the trivialization theorem of Dixmier and Douady \cite{Dix:C*}
of continuous fields with fibers isomorphic to the compact operators. An
important feature of our condition is that it is a priori much weaker than the
condition that  $A$ is $KK_{C(X)}$-equivalent to $C(X)\ot D$. In particular, we
do not need to worry at all about the potentially hard issue of constructing
elements in $KK_{C(X)}(A,C(X)\ot D)$. To illustrate this point, let us note
that it is almost trivial to verify that the local existence of $\sigma$ is
automatic for unital $C(X)$-algebras with fiber $\OOO_n$ and hence to derive
Theorem~\ref{cuntz-algebras-intro}.
 A  C*-algebra $D$ has
 the \emph{automatic local triviality
property}   if any separable  $C(X)$-algebra over
a finite dimensional
 compact
 Hausdorff  space $X$ all of whose fibers are isomorphic
 to $D$ is locally trivial. 
A unital C*-algebra $D$ has
 the \emph{automatic local triviality
property in  the unital sense}  if any separable unital $C(X)$-algebra over
a finite dimensional
 compact
 Hausdorff  space $X$ all of whose fibers are isomorphic
 to $D$ is locally trivial.
The  \emph{automatic triviality property} is defined similarly.
 \begin{theorem}(\emph{Automatic  triviality})\label{o2k}
 A separable continuous $C(X)$-algebra
over a finite dimensional compact Hausdorff space $X$ all of whose fibers
 are isomorphic to $\OOO_2\ot \mathcal{K}$
 is  isomorphic to
 $C(X)\ot \OOO_2\ot \mathcal{K}$.  The C*-algebra $\OOO_2\ot \mathcal{K}$
 is the only
 Kirchberg algebra satisfying the automatic local triviality property and hence the automatic  triviality property.
\end{theorem}
\begin{theorem}\label{A} (\emph{Automatic local triviality in the unital sense})
A unital KK-semiprojective  Kirchberg algebra $D$ has
 the automatic local triviality property in the unital sense if and only if all unital $*$-endomorphisms of
$D$ are KK-equivalences. In that case, if $A$ is a separable unital
$C(X)$-algebra over a finite dimensional
 compact
 Hausdorff  space $X$ all of whose fibers are isomorphic
 to $D$, then $A\cong
C(X)\ot D$ if and only if there is $\sigma \in KK(D,A)$
 such that the induced homomorphism $K_0(\sigma):K_0(D) \to K_0(A)$ maps $[1_D]$ to $[1_A]$.
\end{theorem}
 It is natural to ask if there are other unital Kirchberg algebras besides the
Cuntz algebras which have the automatic local triviality property in the
unital sense. 
Consider the following list $\mathcal{G}$ of pointed abelian groups:

(a) $\big(\{0\},0 \big)$; \quad (b) $(\Z,k)$ with $k>0$;

(c) $\big(\Z/p^{e_1}\oplus \cdots \oplus\Z/p^{e_n},p^{s_1}\oplus \cdots \oplus p^{s_n}\big)$ where
$p$ is a prime, $n\geq 1$, $0\leq s_i < e_i$ for $1\leq i \leq n$ and
$0<s_{i+1}-s_i<e_{i+1}-e_i$ for $1\leq i < n $. If $n=1$ the latter condition is vacuous.
Note that if the integers $1 \leq e_1\leq\cdots \leq e_n$ are given then there exists
integers $s_1,...,s_n$ satisfying the conditions above if and only if $e_{i+1}-e_i\geq 2$
for each $1\leq i \leq n$. If that is the case one can choose $s_i=i-1$ for $1\leq i \leq n$.

(d) $\big(G(p_1)\oplus\cdots \oplus G(p_m), g_1\oplus\cdots \oplus g_m\big)$ where
$p_1,...,p_m$ are distinct primes and each $(G(p_j),g_j)$ is a pointed group
as in (c).

(e)  $\big(\Z\oplus G(p_1)\oplus\cdots \oplus G(p_m), k\oplus g_1\oplus\cdots \oplus g_m\big)$ where  $(G(p_j),g_j)$ are  as in (d).  Moreover we require that $k>0$ is divisible by $p_1^{s_{n(1)}+1}\cdots p_m^{s_{n(m)}+1}$ where $s_{n(j)}$ is defined as in (c) corresponding to the prime $p_j$.

\begin{theorem}(\emph{Automatic local triviality in the unital sense -- the UCT case})
\label{thm-lista-local-triviala} Let $D$ be a unital
Kirchberg algebra which satisfies the UCT and has finitely generated K-theory
groups.  (i) $D$ has the automatic  triviality 
property in the unital sense if and only if $D$ is isomorphic to either $\mathcal{O}_2$ or
$\mathcal{O}_\infty$. (ii) $D$ has the automatic local triviality 
property in the unital sense if and only if $K_1(D)=0$ and $(K_0(D), [1_D])$ is isomorphic to 
one of the pointed groups from the list $\mathcal{G}$.
(iii) If $D$ is as in (ii), then
  a separable unital $C(X)$-algebra $A$ over a
finite dimensional compact Hausdorff space $X$ all of whose fibers are
isomorphic
 to  $D$ is trivial if and only if there exists a  homomorphism of  groups
 $K_0(D)\to K_0(A)$ which maps $[1_D]$ to $[1_A]$.
\end{theorem}

  We use semiprojectivity (in various flavors)
to approximate and represent continuous $C(X)$-algebras as inductive limits of
fibered products of $n$ locally trivial $C(X)$-subalgebras where $n\leq
\mathrm{dim}(X)<\infty$. This  clarifies the local structure of many
$C(X)$-algebras (see Theorem~\ref{basic-approx-K1=torsion-free})
 and gives a new understanding of the K-theory
  of separable  continuous $C(X)$-algebras with arbitrary nuclear fibers.

A  remarkable isomorphism result for  separable
 nuclear strongly purely infinite
stable C*-algebras
 was announced (with an outline of the proof) by
Kirchberg in \cite{Kir:Michael}: two such C*-algebras $A$ and $B$ with the same
primitive spectrum $X$ are isomorphic if and only if they are
$KK_{C(X)}$-equivalent. This is always the case after tensoring with $\OOO_2$.
 However the problem of recognizing when $A$ and $B$ are
$KK_{C(X)}$-equivalent is  open even for very simple spaces $X$ such as the unit
interval or non-Hausdorff spaces with more than two points.

The proof  of Theorem~\ref{basic-approx} (one of our main results) generalizes
and refines a technique that was pioneered for fields over
  zero dimensional spaces in joint work with Pasnicu
   \cite{DadPas:fields-over-zero-dim}
and for fields over an interval in joint work with G. Elliott
\cite{DadEll:bundles-interval}.
 We shall rely heavily
 on the classification theorem (and related results) of
 Kirchberg \cite{Kir:class}
  and Phillips \cite{Phi:class},
 and on the work on non-simple nuclear purely infinite
C*-algebras of  Blanchard and Kirchberg \cite{BK}, \cite{BK:bundles} and
  Kirchberg and R{\o}rdam \cite{KR1}, \cite{KR2}.

The  author is grateful to   E. Blanchard, L.~G.~Brown and N.~C.~Phillips for
useful discussions and comments.

\section{$C(X)$-algebras}
 Let $X$ be a locally
compact Hausdorff space.  A $C(X)$-algebra is a C*-algebra $A$ endowed with a
$*$-\ho\ $\theta$ from $C_0(X)$ to the center $ZM(A)$ of the multiplier algebra
$M(A)$ of $A$ such that $C_0(X)A$ is dense in $A$; see \cite{Kas:inv},
\cite{Blanchard:Hopf}. We write $fa$ rather than $\theta(f)a$ for $f\in C_0(X)$
and $a\in A$.   If $Y\subseteq X$ is a closed set, we let $C_0(X,Y)$ denote the
ideal of $C_0(X)$ consisting of functions vanishing on $Y$. Then $C_0(X,Y)A$ is
a closed two-sided ideal of $A$ (by Cohen factorization). The quotient of $A$
by this ideal is a $C(X)$-algebra denoted by $A(Y)$ and is called the
restriction of $A=A(X)$ to $Y$. The quotient map is denoted by $\pi_Y:A(X)\to
A(Y)$. If $Z$ is a closed subset of $Y$ we have a natural restriction map
$\pi^Y_Z:A(Y)\to A(Z)$ and $\pi_Z=\pi^Y_Z\circ \pi_Y$. If $Y$ reduces to a
point $x$, we write $A(x)$ for $A(\{x\})$ and $\pi_x$ for $\pi_{\{x\}}$. The
C*-algebra $A(x)$ is called the fiber of $A$ at $x$. The image $\pi_x(a)\in
A(x)$ of $a \in A$ is denoted by $a(x)$. A morphism of $C(X)$-algebras $\eta:A
\to B$   induces a morphism $\eta_Y:A(Y)\to B(Y)$. If $A(x)\neq 0$ for  $x$ in a dense subset
of $X$, then $\theta$ is injective. If $X$ is compact, then
$\theta(1)=1_{M(A)}$.
 Let
$A$ be a C*-algebra, $a \in A$ and $\fset,\gset \subseteq A$. Throughout the
paper we will assume that $X$ is a compact Hausdorff space unless stated
otherwise.
If $\ep>0$, we write $a\in_\ep \fset$ if there is $b \in \fset$ such
that $\|a-b\|<\ep$. Similarly, we write $\fset\subset_\ep \gset$ if $a\in_\ep
\gset$ for every $a \in \fset$.  The following lemma collects some basic
properties of $C(X)$-algebras.

\begin{lemma}\label{a(x)delta}
  Let $A$ be a
    $C(X)$-algebra and let $B \subset A$ be a $C(X)$-subalgebra.
    Let $a \in A$ and let $Y$ be a closed  subset of $X$.
    \begin{itemize}
\item[(i)] The map $x\mapsto \|a(x)\|$ is
    upper semi-continuous.
\item[(ii)] $\|\pi_Y(a)\|=\max\{\|\pi_x(a)\|:x\in Y\}$
        \item[(iii)] If $a(x)\in \pi_x(B)$ for all $x \in X$, then $a \in B$.
        \item[(iv)] If $\delta>0$ and $a(x)\in_\delta \pi_x(B)$
        for all $x \in X$, then $a \in_\delta B$.
\item[(v)] The restriction of $\pi_x:A \to A(x)$ to $B$ induces an isomorphism
$B(x)\cong\pi_x(B)$ for all $x \in X$.
    \end{itemize}
\end{lemma}
\begin{proof}
(i), (ii)  are  proved in \cite{Blanchard:Hopf} and (iii) follows from (iv).
(iv): By assumption, for each $x \in X$, there is $b_x \in B$ such that
$\|\pi_x(a-b_x)\|<\delta$. Using (i)  and (ii), we find a closed neighborhood
$U_x$ of $x$ such that $\|\pi_{U_x}(a-b_x)\|<\delta$. Since $X$ is compact,
there is a finite subcover $(U_{x_i})$. Let $(\alpha_i)$ be a partition of
unity subordinated to this cover. Setting $b=\sum_i\,\alpha_i\,b_{x_i}\in B$,
one checks immediately that $\|\pi_x(a-b)\|\leq \sum_i
\alpha_i(x)\|\pi_{x}(a-b_{x_i})\|<\delta$, for all $x \in X$. Thus
$\|a-b\|<\delta$ by (ii). (v): If $\iota:B \hookrightarrow A$ is the inclusion
map, then $\pi_x(B)$ coincides with the image of $\iota_x:B/C(X,x)B \to
A/C(X,x)A$. Thus it suffices to check that $\iota_x$ is injective. If
$\iota_x(b+C(X,x)B)=\pi_x(b)=0$ for some $b\in B$, then $b=fa$ for some $f\in
C(X,x)$ and some $a\in A$. If $(f_\lambda)$ is an approximate unit of $C(X,x)$,
then $b=\lim_\lambda f_\lambda fa =\lim_\lambda f_\lambda b$ and hence $b\in
C(X,x)B$.
 \noindent\end{proof} A $C(X)$-algebra
such that
 the map $x\mapsto \|a(x)\|$  is continuous for all $a \in A$
 is called a \emph{continuous} $C(X)$-algebra or a
C*-bundle \cite{Blanchard:Hopf}, \cite{KirWas:bundles}, \cite{BK:bundles}.
    A
C*-algebra $A$ is a continuous $C(X)$-algebra if and only if $A$ is the
C*-algebra of continuous sections of a continuous field of C*-algebras over $X$
in the sense of \cite[Def. 10.3.1]{Dix:C*},
 (see \cite{Blanchard:Hopf}, \cite{BK:bundles},
\cite{May:cx}).
\begin{lemma}\label{X-is-2nd}
   Let $A$ be a separable continuous  $C(X)$-algebra
   over a locally compact Hausdorff space $X$. If all the fibers of $A$
   are nonzero, then $X$ has a countable basis of open
sets. Thus the compact subspaces of $X$ are metrizable.
\end{lemma}
\begin{proof} Since $A$ is separable, its primitive spectrum $\mathrm{Prim}(A)$
has a countable basis of open sets by \cite[3.3.4]{Dix:C*}. The continuous map
$\eta:\mathrm{Prim}(A)\to X$ (induced by $\theta:C_0(X)\to ZM(A)\cong
C_b(\mathrm{Prim}(A))$) is open since the $C(X)$-algebra $A$ is continuous and
surjective since $A(x)\neq 0$ for all $x\in X$ (see \cite[p.~388]{BK:bundles}
and \cite[Prop.~2.1, Thm.~2.3]{May:cx}).
\end{proof}
\begin{lemma}\label{lemma-auto-continuity} Let $X$ be a compact metrizable
  space.
  A $C(X)$-algebra $A$ all of whose fibers are nonzero and simple
    is  continuous  if and only if there is $e\in A$
    such that $\|e(x)\|\geq 1$ for all $x\in X$.
\end{lemma}
\begin{proof} By Lemma~\ref{a(x)delta}(i) it suffices to
prove that $\liminf_{n\to \infty}\|a(x_n)\|\geq\|a(x_0)\|$ for any $a\in A$ and
any sequence $(x_n)$ converging to $x_0$ in $X$.
  Set $D=A(x_0)$ and
let $e$ be as in the statement.
 Let $\psi:D\to A$ be a set-theoretical  lifting of $\id_D$
 such that $\|\psi(d)\|= \|d\|$ for all $d\in D$. Then $\lim_{n\to \infty}\|\pi_{x_n}\psi(a(x_0))-a(x_n)\|=0$ for
all $a\in A$, by Lemma~\ref{a(x)delta}(i). By applying this to $e$, since
$\|e(x_n)\|\geq 1$, we see that $\liminf_{n\to
\infty}\|\pi_{x_n}\psi(e(x_0))\|\geq 1$.
 Since $D$ is a simple C*-algebra,
if $\varphi_n:D\to B_n$ is a sequence of contractive maps such that $\lim_{n\to
\infty}\|\varphi_n(\lambda c+d)-\lambda\varphi_n(c)-\varphi_n(d)\|=0$,
$\lim_{n\to \infty}\|\varphi_n(cd)-\varphi_n(c)\varphi_n(d)\|=0$, $\lim_{n\to
\infty}\|\varphi_n(c^*)-\varphi_n(c)^*\|=0$,
 for all $c,d\in D$,
$\lambda \in \CCC$, and $\liminf_{n\to \infty}\|\varphi_n(c)\|>0$ for some
$c\in D$, then $\lim_{n\to \infty}\|\varphi_n(c)\|=\|c\|$ for all $c\in D$. In
particular this observation applies to $\varphi_n=\pi_{x_n}\psi$
 by Lemma~\ref{a(x)delta}(i). Therefore
\[\liminf_{n\to
\infty}\|a(x_n)\|\geq \liminf_{n\to
\infty}\big(\|\pi_{x_n}\psi(a(x_0))\|-\|\pi_{x_n}\psi(a(x_0))-a(x_n)\|\big)=\|a(x_0)\|.
\]
Conversely, if $A$ is  continuous, take $e$  to be a large multiple of some
full element of $A$.
\end{proof}
Let $\eta:B \to A$ and $\psi:E \to A$ be $*$-\hos. The \emph{pullback} of these
maps is
\[B\oplus_{\eta,\psi}E=\{(b,e)\in B \oplus E:\,
\eta(b)=\psi(e)\}.\] We are going to use pullbacks in the context of
$C(X)$-algebras. Let $X$ be a compact  space and
  let  $Y,Z $ be closed
subsets of $X$ such that $X=Y\cup Z$.
 The following result is proved in
\cite[Prop.~10.1.13]{Dix:C*} for continuous $C(X)$-algebras.
 \begin{lemma}\label{pullie}
If $A$ is  a
    $C(X)$-algebra, then $A$ is isomorphic to $A(Y)\oplus_{\pi,\pi}A(Z)$, the pullback of the
restriction maps $\pi^Y_{Y\cap Z}:A(Y)\to A(Y\cap Z)$ and $\pi^Z_{Y\cap
Z}:A(Z)\to A(Y\cap Z)$.
\end{lemma}
\begin{proof} By the universal property of pullbacks, the maps $\pi_Y$ and
$\pi_Z$ induce a map $\eta:A \to A(Y)\oplus_{\pi,\pi}A(Z)$,
$\eta(a)=(\pi_Y(a),\pi_Z(a))$, which is injective by Lemma~\ref{a(x)delta}(ii).
Thus it suffices to show that the range of $\eta$ is dense. Let $b,c\in A$ be
such that $\pi_{Y\cap Z}(b-c)=0$ and let $\ep>0$. We shall find $a\in A$ such
that $\|\eta(a)-(\pi_Y(b),\pi_Z(c))\|<\ep$. By Lemma~\ref{a(x)delta}(i), there
is an open neighborhood $V$ of $Y\cap Z$ such that $\|\pi_x(b-c)\|<\ep$ for all
$x\in V$. Let $\{\lambda,\mu\}$ be a partition of unity on $X$ subordinated to
the open cover $\{Y\cup V,Z \cup V\}$. Then $a=\lambda b+\mu c$ is an element
of $A$ which has the desired property.
\end{proof}
 Let $B\subset A(Y)$ and $E\subset A(Z)$ be
$C(X)$-subalgebras such that $ \pi^Z_{Y\cap Z} (E)\subseteq \pi^Y_{Y\cap
Z}(B)$. As an immediate consequence of Lemma~\ref{pullie} we see that the
pullback $B\oplus_{\pi^Z_{Y\cap Z}, \pi^Y_{Y\cap Z}}E$ is isomorphic to the
$C(X)$-subalgebra $B\oplus_{Y\cap Z}E$ of $A$ defined as
\[B\oplus_{Y\cap Z}E=\{a\in A: \pi_Y(a)\in B, \pi_Z(a)\in E\}.\]
\begin{lemma}\label{sumover-Y} The fibers of $B\oplus_{Y\cap Z}E$ are given by
\[\pi_x(B\oplus_{Y\cap Z}E)=\left\{%
\begin{array}{ll}
    {\pi_x(B)}, & \hbox{if $x\in X\setminus Z,$} \\
    {\pi_x(E)}, & \hbox{if $x\in Z,$} \\
\end{array}%
\right.\] and there is an exact sequence of C*-algebras
\begin{equation}\label{use}
    \xymatrix{
 {0}\ar[r]& {\{b\in B:\pi_{Y\cap Z}(b)=0\}}\ar[r]
                & {B\oplus_{Y\cap Z}E} \ar[r]^{\,\,\quad\pi_Z}&{E}\ar[r]&0
}
\end{equation}
\end{lemma}
\begin{proof} Let $x\in X\setminus Z$.  The inclusion
$\pi_x(B\oplus_{Y\cap Z}E)\subset \pi_x(B)$ is obvious by definition. Given
$b\in B$, let us choose $f \in C(X)$ vanishing on $Z$ and such that $f(x)=1$.
Then $a=(fb,0)$ is an element of $A$ by Lemma~\ref{pullie}. Moreover $a\in
B\oplus_{Y\cap Z}E$ and $\pi_x(a)=\pi_x(b)$. We have $\pi_Z(B\oplus_{Y\cap
Z}E)\subset E$, by definition. Conversely, given $e\in E$, let us observe that
$\pi^Z_{Y\cap Z}(e)\in\pi^Y_{Y\cap Z}(B)$ (by assumption) and hence
$\pi^Z_{Y\cap Z}(e)=\pi^Y_{Y\cap Z}(b)$ for some $b\in B$. Then $a=(b,e)$ is an
element of $A$ by Lemma~\ref{pullie} and $\pi_Z(a)=e$. This completes the proof
for the first part of the lemma and also it shows that the map $\pi_Z$ from the
sequence ~\eqref{use} is surjective. Its kernel is identified using
 Lemma~\ref{a(x)delta}(iii).
\end{proof}
Let $X$, $Y$, $Z$ and $A$ be as above. Let $\eta:B\hookrightarrow A(Y)$ be
 a $C(Y)$-linear $*$-monomorphism
     and let $\psi:E \hookrightarrow A(Z)$
     be a $C(Z)$-linear $*$-monomorphism. Assume that
     \begin{equation}\label{inclusion}
   \pi^Z_{Y\cap Z} (\psi(E))\subseteq \pi^Y_{Y\cap Z} (\eta(B)).
\end{equation} This gives a map $\gamma=\eta_{Y\cap Z}^{-1}\psi_{Y\cap Z}: E(Y\cap Z)
\to B(Y\cap Z)$. To simplify notation we let $\pi$ stand for both $\pi^Y_{Y\cap
Z}$ and $\pi^Z_{Y\cap Z}$ in the following lemma.
\begin{lemma}\label{aaa} (a) There
are isomorphisms of $C(X)$-algebras
    \[B\oplus_{\pi,\gamma\pi} E \cong B\oplus_{\pi\eta,\pi\psi}E \cong
     \eta(B)\oplus_{Y\cap Z} \psi (E),\]
    where the second isomorphism is given by the map
    $\chi:B\oplus_{\pi\eta,\pi\psi}E \to A$  induced by the pair
    $(\eta,\psi)$.  Its components $\chi_x$ can be identified with $\psi_x$
    for $x\in Z$ and
    with $\eta_x$ for $x\in X\setminus Z$.

(b) Condition \eqref{inclusion} is equivalent to $\psi(E)\subset
    \pi_Z\big(A\oplus_Y \eta(B)\big)$.

(c) If $\fset$ is a finite subset of $A$ such that $\pi_Y(\fset)\subset_\ep
\eta(B)$ and $\pi_Z(\fset)\subset_\ep \psi(E)$, then $\fset\subset_\ep
\eta(B)\oplus_{Y\cap Z} \psi(E)=\chi\big(B\oplus_{\pi\eta,\pi\psi} E \big)
    $.
\end{lemma}
\begin{proof} This is an immediate corollary of
Lemmas~\ref{a(x)delta},~\ref{pullie},~\ref{sumover-Y}. For illustration, let us
verify (c). By assumption $\pi_x(\fset)\subset_\ep \eta_x(B)$ for all $x \in
X\setminus Z$ and $\pi_z(\fset)\subset_\ep \psi_z(E)$ for all $z \in Z$. We
deduce from Lemma~\ref{sumover-Y} that $\pi_x(\fset)\subset_\ep
\pi_x(\eta(B)\oplus_{Y\cap Z} \psi(E))$
 for all $x \in X$. Therefore $\fset\subset_\ep \eta(B)\oplus_{Y\cap Z}
\psi(E)$ by Lemma~\ref{a(x)delta}(iv).
\end{proof}
\begin{definition}\label{finite-type} Let $\mathcal{C}$ be a class of
 C*-algebras. A  $C(Z)$-algebra $E$ is called
  $\mathcal{C}$-\emph{elementary} if there is a
finite partition of $Z$
    into   closed subsets $Z_1,\dots,Z_r$ ($r\geq
    1$) and there exist C*-algebras $D_1,\dots,D_r$ in $ \mathcal{C}$
     such that $E\cong\bigoplus_{i=1}^r C(Z_i)\ot D_i$.
    The notion of \emph{category} of a  $C(X)$-algebra with respect
   to a
    class $\mathcal{C}$ is defined inductively:
     if $A$ is $\mathcal{C}$-elementary
      then $\mathrm{cat}_\mathcal{C}(A)=0$;
    $\mathrm{cat}_\mathcal{C}(A)\leq n$ if there are closed
    subsets $Y$ and $Z$ of $X$  with
    $X=Y\cup Z$ and there exist  a  $C(Y)$-algebra $B$ such that
     $\mathrm{cat}_\mathcal{C}(B)\leq n-1$,
    a $\mathcal{C}$-elementary
     $C(Z)$-algebra $E$
   and  a  $*$-monomorphism of $C(Y\cap Z)$-algebras
    $\gamma:E(Y\cap Z)\to B(Y\cap
Z)$ such that $A$ is isomorphic to
$$B\oplus_{\,\pi,\gamma\pi_{}} E=\{(b,d)\in B \oplus E:\,
\pi^Y_{Y\cap Z}(b)=\gamma\pi^Z_{Y\cap Z}(d)\}.$$ By definition
$\mathrm{cat}_{\mathcal{C}}(A)=n$ if $n$ is the smallest number with the
property that $\mathrm{cat}_{\mathcal{C}}(A)\leq n.$ If  no such $n$ exists,
then $\mathrm{cat}_{\mathcal{C}}(A)=\infty.$
\end{definition}
\begin{definition}\label{stack}
Let $\mathcal{C}$ be a class of C*-algebras and let $A$ be a $C(X)$-algebra. An
$n$-\emph{fibered $\mathcal{C}$-monomorphism $(\psi_0,\dots,\psi_n)$ into $A$}
      consists of $(n+1)$ $*$-\mos\ of $C(X)$-algebras
$\psi_i:E_i \to A(Y_i)$, where $Y_0,\dots,Y_n$ is a closed cover of $X$,
    each $E_i$ is a
   $\mathcal{C}$-elementary  $C(Y_i)$-algebra and
    \begin{equation}\label{containement-fibered}
\pi^{Y_i}_{Y_i\cap
    Y_j}\psi_i(E_i)\subseteq \pi^{Y_j}_{Y_i\cap
    Y_j}\psi_j(E_j), \quad \text{ for all}\,\, i\leq j.
\end{equation}
    Given an $n$-fibered morphism into $A$ we have an associated \emph{continuous} $C(X)$-algebra
    defined  as the  fibered product (or pullback)  of the $*$-\mos\ $\psi_i$:
    \begin{equation}\label{n-pull}
    A(\psi_0,\dots,\psi_n)=\{(d_0,\dots d_n)\,:\,\, d_i\in E_i,\,
    \pi^{Y_i}_{Y_i\cap Y_j}\psi_i(d_i)=\pi^{Y_j}_{Y_i\cap Y_j}\psi_j(d_j)
    \,\,\text{for all}\, i,j\}
\end{equation}
and an induced $C(X)$-\mo\ (defined by using  Lemma~\ref{pullie})
$$\eta=\eta_{(\psi_0,\dots,\psi_n)}:A(\psi_0,\dots,\psi_n)\to A\subset
\bigoplus_{i=0}^n A(Y_i), $$
\[\eta(d_0,\dots d_n)=\big(\psi_0(d_0),\dots,\psi_n(d_n)\big).\]
 There are natural
coordinate maps $p_i:A(\psi_0,\dots,\psi_n)\to E_i$, $p_i(d_0,\dots,d_n)=d_i$.
 Let us set
$X_k=Y_k\cup\cdots\cup Y_n$. Then $(\psi_k,\dots\psi_n)$ is an $(n-k)$-fibered
$\mathcal{C}$-monomorphism into $A(X_k)$. Let $\eta_k :
A(X_k)(\psi_k,\dots\psi_n)\to A(X_k)$ be the induced map and set
$B_k=A(X_k)(\psi_k,\dots\psi_n)$. Let us note that $B_0=A(\psi_0,\dots,\psi_n)$
and that there are natural $C(X_{k-1})$-isomorphisms
\begin{equation}\label{eq-recurrence-B(k)}
B_{k-1}\cong B_k\oplus_{\pi\eta_k,\pi\psi_{k-1}} E_{k-1}\cong
B_k\oplus_{\pi,\gamma_k\pi} E_{k-1}.
\end{equation}
where $\pi$ stands for $\pi_{X_k\cap Y_{k-1}}$ and $\gamma_k: E_{k-1}(X_k\cap
Y_{k-1}) \to B_k(X_k\cap Y_{k-1})$ is defined by
$(\gamma_k)_x=(\eta_k)_x^{-1}(\psi_{k-1})_x$, for all $x\in X_k\cap Y_{k-1}$.
In particular, this decomposition shows that
$\mathrm{cat}_{\mathcal{C}}(A(\psi_0,\dots,\psi_n))\leq n$.
\end{definition}
\begin{lemma}\label{lemma-needed-properties-for-B}
Suppose that the class $\mathcal{C}$ from Definition~\ref{finite-type}
 consists
of stable Kirchberg algebras. If $A$ is a $C(X)$-algebra over a compact
 metrizable space $X$ such that
$\mathrm{cat}_{\mathcal{C}}(A)<\infty$, then $A$ contains a full properly
infinite projection and $A\cong A\ot\OOO_\infty\ot \mathcal{K}$.
\end{lemma}
\begin{proof} We prove this by induction on $n=\mathrm{cat}_{\mathcal{C}}(A)$.
The case $n=0$ is immediate since $D\cong D\ot \OOO_\infty$ for any Kirchberg
algebra $D$  \cite{Kir:class}. Let $A=B\oplus_{\pi,\gamma\pi} E$ where $B$, $E$
and $\gamma$ are as in Definition~\ref{finite-type} with
$\mathrm{cat}_{\mathcal{C}}(B)= n-1$ and $\mathrm{cat}_{\mathcal{C}}(E)= 0$.
Let us consider the exact sequence $0\to J\to A \to E \to 0$, where $J=\{b\in
B:\pi_{Y\cap Z}(b)=0\}$.
 Since $J$ is an ideal of $B\cong B\ot \OOO_\infty\ot \mathcal{K}$,
 $J$ absorbs $\OOO_\infty\ot \mathcal{K}$ by
 \cite[Prop.~8.5]{KR2}. Since both $E$ and
 $J$ are stable and purely infinite, it follows that $A$ is stable
 by \cite[Prop.~6.12]{Ror:stable-Japan} and purely infinite by
  \cite[Prop.~3.5]{KR2}.  Since $A$
 has Hausdorff primitive spectrum, $A$ is strongly purely infinite by
\cite[Thm.~5.8]{BK}. It follows that
 $A\cong A\ot \OOO_\infty$ by
\cite[Thm.~9.1]{KR2}.
  Finally
 $A$ contains a full properly infinite projection
since there is  a full embedding of $\OOO_2$  into $A$ by \cite[Prop.~5.6]{BK}.
\end{proof}

\section{Semiprojectivity}\label{section:Semiprojectivity}
In this section we study the notion of $KK$-semiprojectivity. The main result
is Theorem~\ref{equivalences}.
 Let
$A$ and $B$ be C*-algebras. Two $*$-\hos\ $\varphi,\psi:A \to B$ are
approximately unitarily equivalent, written $\varphi\approx_u \psi,$ if  there
is a sequence of unitaries $(u_n)$ in the C*-algebra $B^+=B+\CCC 1$ obtained by
adjoining a unit to $B$, such that $\lim_{n \to
\infty}\|u_n\varphi(a)u_n^*-\psi(a)\|=0$ for all $a \in A$. We say that
 $\varphi$ and $\psi$ are asymptotically unitarily equivalent,
 written $\varphi\approx_{uh}\psi,$ if there is a norm continuous
 unitary valued
 map $t \to u_t\in B^+$, $t\in [0,1)$, such that
 $\lim_{t \to 1}\|u_t\varphi(a)u_t^*-\psi(a)\|=0$
 for all $a \in A$.
 A $*$-\ho\ $\varphi:D\to A$ is full if $\varphi(d)$  is not
contained in any proper two-sided closed ideal of $A$ if $d\in D$ is nonzero.

 We shall use several times
  Kirchberg's Theorem \cite[Thm.~8.3.3]{Ror:encyclopedia} and  the
following theorem of Phillips~\cite{Phi:class}.
\begin{theorem}\label{Phillips-basic}
    Let $A$ and $B$ be  separable C*-algebras such that $A$ is simple and
    nuclear,
     $B\cong B\ot \OOO_\infty$, and there
     exist full projections $p\in A$ and $q\in B$. For any
    $\sigma\in KK(A, B)$
    there is a full $*$-\ho\ $\varphi:A\to B$ such that $KK(\varphi)=\sigma$.
    If $K_0(\sigma)[p]=[q]$ then we may arrange that $\varphi(p)=q$.
    If $\psi:A \to B$ is another  $*$-\ho\ such that
    $KK(\psi)=KK(\varphi)$ and $\psi(p)=q$, then $\varphi\approx_{uh}\psi$
     via a path of
    unitaries
    $t\mapsto u_t\in U(qBq)$.
\end{theorem}
 Theorem~\ref{Phillips-basic} does not appear in this form in
 \cite{Phi:class} but it is an immediate consequence of
\cite[Thm.~4.1.1]{Phi:class}. Since $pAp\ot\mathcal{K}\cong A\ot \mathcal{K}$
and $qBq\ot\mathcal{K}\cong B\ot \mathcal{K}$ by \cite{Bro:stabher}, and
$qBq\ot \OOO_\infty\cong qBq$ by \cite[Prop.~8.5]{KR2}, it suffices to discuss
the case when $p$ and $q$ are the units of $A$ and $B$. If $\sigma$ is given,
\cite[Thm.~4.1.1]{Phi:class} yields a full $*$-\ho\ $\varphi:A \to B\ot
\mathcal{K}$ such that $KK(\varphi)=\sigma$. Let $e\in \mathcal{K}$ be a
rank-one projection and suppose that $[\varphi(1_A)]=[1_B\ot e]$ in $K_0(B)$.
Since both $\varphi(1_A)$ and $1_B\ot e$ are full projections and $B\cong B\ot
\OOO_\infty$, it follows by \cite[Lemma~2.1.8]{Phi:class} that
$u\varphi(1_A)u^*=1_B\ot e$ for some unitary in $(B\ot \mathcal{K})^+$.
Replacing $\varphi$ by $u\,\varphi\, u^*$ we can arrange that
$KK(\varphi)=\sigma$ and $\varphi$ is unital. For the second part of the
theorem let us note that
 any unital $*$-\ho\ $\varphi:A\to B$ is full and  if two unital $*$-\hos\
$\varphi,\psi:A\to B$  are asymptotically unitarily equivalent when regarded as
maps into $B\otimes \mathcal{K}$, then  $\varphi\approx_{uh}\psi$ when regarded
as maps into $B$, by an argument from the proof of
\cite[Thm.~4.1.4]{Phi:class}.

A separable nonzero C*-algebra $D$ is \emph{semiprojective}  \cite{Bla:shape}
if for any separable C*-algebra $A$ and any increasing sequence of two-sided
closed ideals $(J_n)$ of $A$ with $J=\overline{\cup_n J_n}$, the natural map
$\varinjlim\, \mathrm{Hom}(D,A/J_n)\to \mathrm{Hom}(D,A/J)$ (induced by
$\pi_n:A/J_n\to A/J$) is surjective.
 If we weaken this condition and require only that the above map has dense range,
  where $\mathrm{Hom}(D,A/J)$ is given
  the point-norm topology, then $D$ is called
\emph{weakly semiprojective} \cite{Eil-Lor:contingencies}. These definitions do
not change if we drop the separability of $A$.
 We shall use (weak)
semiprojectivity in the following context. Let $A$ be a $C(X)$-algebra (with
$X$ metrizable), let $x\in X$ and set $U_n=\{y\in X:\, d(y,x)\leq 1/n\}$. Then
$J_n=C(X,U_n)A$ is an increasing sequence of ideals of $A$ such that
$J=C(X,x)A$, $A/J_n\cong A(U_n)$ and $A/J\cong A(x)$.
\begin{examples}\label{Class-C} (Weakly semiprojective C*-algebras)
  Any finite dimensional C*-algebra is semiprojective.
   A  Kirchberg
   algebra $D$ satisfying the UCT and
   having finitely generated K-theory groups is \wsp\
   by work of  Neub\"user \cite{Neu:thesis}, H. Lin \cite{Lin:wsp} and Spielberg
   \cite{spielberg:wsp}. This also follows from Theorem~\ref{equivalences}
    and Proposition~\ref{fg-uct} below.
If in addition $K_1(D)$ is torsion free, then $D$ is semiprojective
      as proved by Spielberg
   \cite{Spi:semiproj} who extended the foundational
  work of Blackadar \cite{Bla:shape} and Szymanski \cite{Szy:semiproj}.
\end{examples}
The following generalizations of two results of Loring \cite{Lor:lifting} are
used in section~\ref{three}; see \cite{DadEll:bundles-interval}.
\begin{proposition}\label{approx-lift>exact-lift}
    Let $D$ be a separable semiprojective C*-algebra.
    For any finite set $\fset \subset
D$ and any $\ep>0$, there exist a finite set $\gset \subset D$ and $\delta>0$
with the following property.  Let $\pi:A \to B$ be a surjective $*$-\ho, and
let $\varphi:D \to B$ and $\gamma:D \to A$ be $*$-\hos\ such that
$\|\pi\gamma(d)-\varphi(d)\|<\delta$ for all $d \in \gset$. Then there is a
$*$-\ho\ $\psi:D \to A$ such that $\pi\psi=\varphi$ and
$\|\gamma(c)-\psi(c)\|<\ep$ for all $c \in \fset$.
\end{proposition}
\begin{proposition}\label{homotopy-semiprojective}
Let $D$ be a separable semiprojective C*-algebra.
    For any finite set $\fset \subset
D$ and any $\ep>0$, there exist a finite set $\gset \subset D$ and $\delta>0$
with the following property. For any two $*$-\hos\ $\varphi,\psi:D \to B$ such
that $\|\varphi(d)-\psi(d)\|<\delta$ for all $d\in \gset$, there is a homotopy
$\Phi\in \mathrm{Hom}(D, C[0,1]\ot B)$ such that $\Phi_0=\varphi$ to
$\Phi_1=\psi$ and $\|\varphi(c)-\Phi_t(c)\|<\ep$ for all $c\in\fset$ and $t\in
[0,1]$.
\end{proposition}

\begin{definition} \label{KK-stable-true} A separable C*-algebra $D$ is
 $KK$-stable if  there is a finite set
   $\gset \subset D$ and there is $\delta>0$ with the property that
  for any two  $*$-\hos\ $\varphi,\psi:D \to A$ such that
   $\|\varphi(a)-\psi(a)\|<\delta$ for all $a\in \gset$, one has
   $KK(\varphi)=KK(\psi)$.
\end{definition}

\begin{corollary}\label{SP-implies-KK-stable}
Any  semiprojective  C*-algebra is \wsp\ and $KK$-stable.
\end{corollary}
\begin{proof} This follows from Proposition~\ref{homotopy-semiprojective}.
\end{proof}
\begin{proposition}\label{char-semiprojectivity}
Let $D$ be a separable \wsp\ C*-algebra.  For any finite set $\fset \subset D$
and any $\ep>0$ there exist a finite set $\gset \subset D$ and $\delta>0$ such
that for any C*-algebras $B\subset A$ and any $*$-\ho\ $\varphi:D \to A$ with
$\varphi(\gset)\subset_\delta B$, there is a $*$-\ho\ $\psi:D \to B$ such that
$\|\varphi(c)-\psi(c)\|<\ep$ for all $c \in \fset$. If in addition $D$ is
 $KK$-stable, then we can choose $\gset$ and
$\delta$ such that we also have $KK(\psi)=KK(\varphi)$.
\end{proposition}
\begin{proof} This follows from \cite[Thms.~3.1,
4.6]{Eil-Lor:contingencies}. Since the result is essential to us we include a
short proof. Fix $\fset$ and $\ep$.
 Let
$(\gset_n)$ be an increasing sequence  of finite subsets of $D$ whose union is
dense in $D$. If the statement is not true, then there are sequences of
C*-algebras $C_n\subset A_n$ and
 $*$-\hos\ $\varphi_n:D\to A_n$ satisfying
 $\varphi_n(\gset_n)\subset_{1/n} C_n$ and with the property that
 for any  $n\geq 1$ there is no $*$-\ho\ $\psi_n:D\to C_n$
such that $\|\varphi_n(c)-\psi_n(c)\|<\ep$
  for all $c\in \fset$. Set $B_i=\prod_{n\geq i} A_n$ and $E_i=\prod_{n\geq i}
  C_n\subset B_i$.
If $\nu_i:B_i\to B_{i+1}$ is the natural projection, then $\nu_i (E_i)=
E_{i+1}.$ Let us observe that if we define
 $\Phi_i:D\to B_i$
by $\Phi_i(d)=(\varphi_i(d),\varphi_{i+1}(d),\dots)$, then the image of
$\Phi=\varinjlim\,\Phi_i:D\to \varinjlim\,(B_i,\nu_i)$
 is contained in $\varinjlim\,(E_i,\nu_i)$.
  Since $D$ is weakly semiprojective, there is $i$ and a $*$-\ho\
  $\Psi_i:D \to E_i$, of the form $\Psi_i(d)=(\psi_i(d),\psi_{i+1}(d),\dots)$
such that $\|\Phi_i(c)-\Psi_i(c)\|<\ep$ for all $c\in \fset$. Therefore
$\|\varphi_i(c)-\psi_i(c)\|<\ep$ for all $c\in \fset$ which gives a
contradiction.
\end{proof}
It is useful to combine Propositions~\ref{char-semiprojectivity} and ~\ref{approx-lift>exact-lift} in a single statement.
\begin{proposition}\label{semiproj_combined}
  Let $D$ be a separable semiprojective C*-algebra.
    For any finite set $\fset \subset
D$ and any $\ep>0$, there exist a finite set $\gset \subset D$ and $\delta>0$
with the following property.  Let $\pi:A \to B$ be a surjective $*$-\ho\ which maps
a C*-subalgebra $A'$ of $A$ onto a C*-subalgebra $B'$ of $B$. 
Let $\varphi:D \to B'$ and $\gamma:D \to A$ be $*$-\hos\ such that $\gamma(\gset)\subset_{\delta} A'$ and 
$\|\pi\gamma(d)-\varphi(d)\|<\delta$ for all $d \in \gset$. Then there is a
$*$-\ho\ $\psi:D \to A'$ such that $\pi\psi=\varphi$ and
$\|\gamma(c)-\psi(c)\|<\ep$ for all $c \in \fset$.
\end{proposition}
\begin{proof}
 Let $\gset_L$ and $\delta_L$ be given by Proposition~\ref{approx-lift>exact-lift} applied
to the input data $\fset$ and $\ep/2$. We may assume that $\fset \subset \gset_L$ and
$\ep>\delta_L$. Next, let $\gset_P$ and $\delta_P$ be given by Proposition~\ref{char-semiprojectivity} applied
to the input data $\gset_L$ and $\delta_L/2$. We show now that $\gset:=\gset_L \cup \gset_P$
and $\delta:=\min\{\delta_P,\delta_L/2\}$ have the desired properties.
We have $\gamma(\gset_P)\subset_{\delta_P} A'$ since $\gset_P\subset \gset$ and $\delta\leq \delta_P$. By Proposition~\ref{char-semiprojectivity} there is a $*$-\ho\ $\gamma':D\to A'$
such that $\|\gamma'(d)-\gamma(d)\|<\delta_L/2$ for all $d\in \gset_L$. Then,
since $\gset_L\subset \gset$ and $\delta\leq \delta_L/2$,
\[\|\pi\gamma'(d)-\varphi(d)\|\leq \|\pi\gamma'(d)-\pi\gamma(d)\|+\|\pi\gamma(d)-\varphi(d)\|
<\delta_L/2+\delta \leq \delta_L\]
for all $d\in \gset_L$.
Therefore we can invoke Proposition~\ref{approx-lift>exact-lift} to perturb $\gamma'$
to a $*$-\ho\ $\psi:D \to A'$ such that $\pi \psi=\varphi$ and $\|\gamma'(d)-\psi(d)\|<\ep/2$
for all $d\in \fset$. Finally we observe that for $d\in \fset \subset \gset_L$
\[\|\gamma(d)-\psi(d)\|\leq \|\gamma(d)-\gamma'(d)\|+\|\gamma'(d)-\psi(d)\|<\delta_L/2+\ep/2<\ep.\]
\end{proof}

\begin{definition} \label{KK-sp+continuous}
(a) A separable  C*-algebra $D$ is \emph{KK-semiprojective}  if for any
separable C*-algebra $A$ and any increasing sequence of two-sided closed ideals
$(J_n)$ of $A$ with $J=\overline{\cup_n J_n}$, the natural map $\varinjlim\,
KK(D,A/J_n)\to KK(D,A/J)$ is surjective.

(b) We say that the functor $KK(D,-)$ is \emph{continuous} if for any inductive
system $B_1\to B_2\to...$  of separable  C*-algebras, the induced map
 $\varinjlim\, KK(D,B_n)\to KK(D,\varinjlim\,B_n)$
 is bijective.
\end{definition}

\begin{proposition}\label{ASP-implies-KK-stable}
   Any separable KK-semiprojective C*-algebra  is
   $KK$-stable.
\end{proposition}
\begin{proof} We shall prove the statement by contradiction. Let $D$ be
separable KK-semiprojective C*-algebra. Let $(\gset_n)$ be an increasing
sequence  of finite subsets of $D$ whose union is dense in $D$. If the
statement is not true, then there are sequences of
 $*$-\hos\ $\varphi_n,\psi_n:D\to A_n$ such that
$\|\varphi_n(d)-\psi_n(d)\|<1/n$ for all $d\in \gset_n$ and yet
$KK(\varphi_n)\neq KK(\psi_n)$ for all $n\geq 1$. Set $B_i=\prod_{n\geq i} A_n$
and let $\nu_i:B_i\to B_{i+1}$ be the natural projection. Let us define
$\Phi_i,\Psi_i:D\to B_i$ by $\Phi_i(d)=(\varphi_i(d),\varphi_{i+1}(d),\dots)$
and $\Psi_i(d)=(\psi_i(d),\psi_{i+1}(d),\dots)$, for all $d$ in $D$. Let $B_i'$
be the separable C*-subalgebra of $B_i$ generated by the images of $\Phi_i$ and
$\Psi_i$. Then $\nu_i(B'_i)= B'_{i+1}$ and
 one
verifies immediately that $\varinjlim\,\Phi_i=\varinjlim\,\Psi_i:D\to
\varinjlim\,(B'_i,\nu_i)$. Since $D$ is $KK$-semiprojective, we must have
$KK(\Phi_i)=KK(\Psi_i)$ for some $i$ and hence $KK(\varphi_n)=KK(\psi_n)$ for
all $n\geq i$. This gives a contradiction.
\end{proof}
\begin{proposition}\label{propAotKwsp} A unital Kirchberg
algebra $D$  is $KK$-stable if and only if
    $D\ot \mathcal{K}$ is  $KK$-stable.  $D$ is \wsp\  if and only if
    $D\ot \mathcal{K}$ is \wsp.
\end{proposition}
\begin{proof} Since $KK(D,A)\cong KK(D,A\ot \mathcal{K})
\cong KK(D\ot \mathcal{K},A\ot \mathcal{K})$ the first part of the proposition
is immediate. Suppose now that  $D\ot \mathcal{K}$ is \wsp.  Then $D$ is \wsp\
as shown in the proof of \cite[Thm.~2.2]{spielberg:wsp}. Conversely, assume
that $D$ is \wsp. It suffices to find
$\alpha\in\mathrm{Hom}(D\ot\mathcal{K},D)$ and a sequence $(\beta_n)$ in
$\mathrm{Hom}(D ,D\ot\mathcal{K})$ such that $\beta_n\alpha$ converges to
$\id_{D\ot\mathcal{K}}$ in the point-norm topology.
 Let $s_i$ be the canonical generators of $\OOO_\infty$. If
$(e_{ij})$ is a system of matrix units for $\mathcal{K}$, then
$\lambda(e_{ij})=s_is^*_j$ defines a  $*$-\ho\ $\mathcal{K}\to \OOO_\infty$
such that $KK(\lambda)\in KK(\mathcal{K},\OOO_\infty)^{-1}$. Therefore, by
composing $\id_{D}\ot \lambda$ with some isomorphism $ D\ot \OOO_\infty\cong D
$ (given by \cite[Thm.~7.6.6]{Ror:encyclopedia}) we obtain a $*$-\mo\
$\alpha:D\ot \mathcal{K}\to D$ which induces a KK-equivalence. Let $\beta :
D\to D\ot \mathcal{K}$ be defined by $\beta(d)= d\ot e_{11}$. Then
$\beta\alpha\in \End(D\ot\mathcal{K})$ induces a KK-equivalence and hence after
replacing $\beta$ by $\theta\beta$ for some automorphism $\theta$ of
$D\ot\mathcal{K}$, we may arrange that $KK(\beta\alpha)=KK(\id_D)$. By
Theorem~\ref{Phillips-basic}, $\beta\alpha\approx_u \id_{D\ot\mathcal{K}},$ so
that there is a sequence of unitaries $u_n\in (D\ot\mathcal{K})^+$ such that
$u_n\beta\alpha(-) u_n^*$ converges to $\id_{D\ot\mathcal{K}}$.\end{proof}
\begin{theorem}\label{equivalences}
    For a separable  C*-algebra $D$ consider the following
    properties:

   (i) $D$ is $KK$-semiprojective.

   (ii) The functor $KK(D,-)$ is continuous.

   (iii) $D$ is \wsp\ and $KK$-stable.

        Then (i) $\Leftrightarrow$ (ii). Moreover, (iii) $\Rightarrow$ (i) if
    $D$ is  nuclear and (i) $\Rightarrow$ (iii) if $D$ is a
    Kirchberg algebra. Thus (i) $\Leftrightarrow$ (ii)  $\Leftrightarrow$ (iii)
    for any
    Kirchberg algebra $D$.
\end{theorem}\begin{proof} The implication (ii) $\Rightarrow$ (i) is obvious.
(i) $\Rightarrow$ (ii): Let $(B_n,\gamma_{n,m})$ be an inductive system with
inductive limit $B$ and let $\gamma_n:B_n\to B$ be the canonical maps. We have
an induced  map $\beta:\varinjlim\, KK(D,B_n)\to KK(D,B)$. First we show that
$\beta$ is surjective.  The mapping telescope construction of L.~G.~Brown (as
described in the proof of \cite[Thm.~3.1]{Bla:shape}) produces an inductive
system of C*-algebras $(T_n,\eta_{n,m})$ with inductive limit $B$ such that
each $\eta_{n,n+1}$ is surjective, and each canonical map $\eta_n:T_n\to B$ is
homotopic to $\gamma_n\alpha_n$ for some $*$-\ho\ $\alpha_n:T_n\to B_n$. In
particular $KK(\eta_n)=KK(\gamma_n)KK(\alpha_n)$.
 Let $x \in
KK(D,B)$. By (i) there are $n$ and $y\in KK(D,T_n)$ such that $KK(\eta_n)y=x$
and hence $KK(\gamma_n)KK(\alpha_n)y=x$. Thus
 $z=KK(\alpha_n)y\in KK(D,B_n)$ is a lifting
of $x$. Let us show now that the map $\beta$ is injective. Let $x$ be an
element in the kernel of the map $ KK(D,B_n)\to KK(D,B).$
 Consider the
commutative diagram whose exact rows are portions of the Puppe sequence in
KK-theory \cite[Thm.~19.4.3]{Bla:k-theory} and with vertical maps induced by
$\gamma_m:B_m\to B$, $m\geq n$.
\[
\xymatrix{
 {KK(D,C_{\gamma_n})}\ar[r]& {KK(D,B_n)}
 \ar[r]& {KK(D,B)} \\
{KK(D,C_{\gamma_{n,m}})}\ar[r]\ar[u]& {KK(D,B_n)}\ar@{=}[u]
 \ar[r]& {KK(D,B_m)}\ar[u]\\
  }\]
By  exactness, $x$ is the image of some element $y\in KK(D,C_{\gamma_n})$.
Since $C_{\gamma_n}=\varinjlim C_{\gamma_{n,m}}$,  the map $\varinjlim
KK(D,C_{\gamma_{n,m}})\to KK(D,C_{\gamma_n})$ is surjective by the first part
of the proof. Therefore there is $m\geq n$ such that $y$ lifts to some
 $z \in KK(D,C_{\gamma_{n,m}})$.
  The image
 of $z$ in $KK(D,B_m)$ equals $KK(\gamma_{n,m})x$ and vanishes
  by exactness of the bottom row.

(iii) $\Rightarrow$ (i): Let $A$, $(J_n)$ and $J$ be as in
Definition~\ref{KK-sp+continuous}.
 Using the five-lemma and the split exact sequence $0\to KK(D,A)\to
KK(D,A^+)\to KK(D,\CCC)\to 0$, we reduce the proof to the case when $A$ is
unital. Let $x\in KK(D,A/J)$. Since the map
   $KK(D^+,A/J)\to KK(D,A/J)$ is surjective, $x$ lifts to some element $x^+\in
   KK(D^+,A/J)$.
 By
   \cite[Thm.~8.3.3]{Ror:encyclopedia}, since $D^+$ is nuclear,
   there is a $*$-\ho\
   $\Phi:D^+\to A/J\ot \OOO_\infty\ot \mathcal{K}$ such that
   $KK(\Phi)=x^+$ and hence if set $\varphi=\Phi|_D$, then $KK(\varphi)=x$.  Since $D$ is \wsp, there are $n$ and a $*$-\ho\
   $\psi:D\to A/J_n\ot \OOO_\infty\ot \mathcal{K}$ such that
   $\|\pi_n\psi(d)-\varphi(d)\|<\delta$ for all $d\in\gset$,
   where $\gset$ and $\delta$ are
   as in the definition of $KK$-stability. Therefore
   $KK(\pi_n\psi)=KK(\varphi)$ and hence $KK(\psi)$ is a lifting of $x$ to
   $KK(D,A/J_n)$.

(i) $\Rightarrow$ (iii): $D$ is $KK$-stable by
Proposition~\ref{ASP-implies-KK-stable}. It remains to show that $D$ is \wsp.
Since any nonunital Kirchberg algebra is isomorphic to the stabilization of a
unital one (see \cite[Prop.~4.1.3]{Ror:encyclopedia}) and since by
Proposition~\ref{propAotKwsp} $D$ is KK-semiprojective if and only if $D\ot
\mathcal{K}$ is KK-semiprojective,  we may assume that $D$ is unital. Let $A$,
$(J_n)$, $\pi_{m,n}:A/J_m\to A/J_n$ $(m \leq n$) and $\pi_n:A/J_n\to A/J$ be as
in the definition of weak semiprojectivity. By \cite[Cor.~2.15]{Bla:shape}, we
may assume that $A$ and the $*$-\ho\ $\varphi:D\to A$ (for which we want to
construct an approximative lifting) are unital. In particular $\varphi$ is
injective since $D$ is simple. Set $B=\varphi(D)\subset A/J$ and
$B_n=\pi_n^{-1}(B)\subset A/J_n$. The corresponding maps $\pi_{m,n}:B_m\to B_n$
$(m \leq n$) and $\pi_n:B_n\to B$ are surjective and they induce an isomorphism
$\varinjlim \,(B_n,\pi_{n,n+1})\cong B$.

Given $\ep>0$ and $\fset\subset D$ (a finite set) we are going to produce an
approximate lifting $\varphi_n:D\to B_n$ for $\varphi$. Since $1_B$ is a
properly infinite projection, it follows by
 \cite[Props.~2.18 and 2.23]{Bla:shape} that the unit
    $1_{n}$ of $B_n$ is a properly infinite projection, for all sufficiently
   large $n$.
    Since $D$ is
$KK$-semiprojective, there exist $ m$ and an element $h\in KK(D,B_m)$ which
lifts $KK(\varphi)$  such that $K_0(h)[1_D]=[1_{m}]$.
 By
   \cite[Thm.~8.3.3]{Ror:encyclopedia}, there is a full $*$-\ho\
   $\eta:D\to B_m\ot \mathcal{K}$ such that $KK(\eta)=h$.
   By \cite[Prop.~4.1.4]{Ror:encyclopedia}, since both $\eta(1_D)$ and $1_{m}$
   are full and properly infinite projections in $B_m\ot \mathcal{K}$,
   there is a partial isometry $w\in B_m\ot \mathcal{K}$ such that
   $w^*w=\eta(1_D)$ and $ww^*=1_{m}$. Replacing $\eta$ by $w\,\eta(-)\, w^*$,
   we may assume that $\eta:D\to B_m$ is unital.
    Then $KK(\pi_m\eta)=KK(\pi_m)h=KK(\varphi)$.
    By Theorem~\ref{Phillips-basic}, $\pi_m\eta\approx_{uh}\varphi$.
    Thus
   there is a unitary $u\in B$
   such that $\|u\pi_m\eta(d)u^*-\varphi(d)\|<\ep$ for all $d\in \fset$.
   Since $C(\mathbb{T})$ is semiprojective, there is $n\geq m$
   such that $u$ lifts to a 
    unitary $u_n\in B_n$.
   Then $\varphi_n:=u_n\,\pi_{m,n}\eta(-)\,u_n^*$ is a $*$-\ho\ from $D$ to $B_n$ such that
$\|\pi_n\varphi_n(d)-\varphi(d)\|<\ep$ for all $d\in \fset$.
   \end{proof}

   \begin{corollary}\label{SP-implies-ASP}
  Any separable nuclear semiprojective  C*-algebra is
   $KK$-semiprojective.
\end{corollary}
\begin{proof}
This is very similar to the proof of the implication (iii) $\Rightarrow$ (i) of
Theorem~\ref{equivalences}. Alternatively, the statement follows from Corollary
~\ref{SP-implies-KK-stable} and Theorem~\ref{equivalences}.
\end{proof}
   Blackadar  has shown that a semiprojective Kirchberg algebra satisfying
   the UCT
   has finitely generated K-theory groups \cite[Prop.~8.4.15]{Ror:encyclopedia}.
   A similar argument gives the following:
   \begin{proposition}\label{fg-uct}
    Let $D$ be a  separable C*-algebra satisfying the UCT.
    Then $D$ is $KK$-semiprojective if and only $K_*(D)$ is finitely
    generated.
\end{proposition} \begin{proof}If $K_*(D)$ is finitely generated, then  $D$
 is $KK$-semiprojective by \cite{RosSho:UCT}.
 Conversely, assume that $D$ is $KK$-semiprojective.
 Since $D$ satisfies the UCT, we infer that if $G=K_i(D)$ ($i=0,1$), then $G$ is
 semiprojective in the category of countable abelian groups, in the sense that
 if $H_1\to H_2\to\cdots$ is an inductive system of countable abelian groups
 with inductive limit $H$, then the natural map
 $\varinjlim\, \mathrm{Hom}(G,H_n)\to \mathrm{Hom}(G,H)$ is surjective.
 This implies that $G$ is finitely generated. Indeed, taking $H=G$, we
 see that $\id_G$ lifts to  $\mathrm{Hom}(G,H_n)$ for some finitely generated
 subgroup $H_n$ of $G$ and hence $G$ is a quotient of $H_n$.
\end{proof}

\section{Approximation of $C(X)$-algebras}

In this section we use weak semiprojectivity to approximate  a continuous $C(X)$-algebra
$A$ by $C(X)$-subalgebras given by pullbacks of $n$-fibered monomorphisms into
$A$.\begin{lemma}\label{fullness}
    Let $D$ be a  finite direct sum of simple C*-algebras
    and let $\varphi,\psi:D \to A$ be
    $*$-\hos. Suppose that $\hset\subset D$ contains a nonzero element from
    each simple direct summand of $D$.  If
    $\|\psi(d)-\varphi(d)\|\leq\|d\|/2$ for all $d\in \hset$,
     then $\varphi$ is injective if and only if $\psi$ is injective.
\end{lemma}
\begin{proof} Let us note that  $\varphi$ is injective if and only if
 $\|\varphi(d)\|=\|d\|$ for  all $d\in\hset$. Therefore if $\varphi$ is
 injective, then
   $\|\psi(d)\|\geq \|\varphi(d)\|-\|\psi(d)-\varphi(d)\|\geq \|d\|/2$
   for all $d\in\hset$
and hence  $\psi$ is injective.\end{proof}
 A sequence $(A_n)$ of
subalgebras of $A$ is called \emph{exhaustive} if for any finite subset $\fset$
of $A$ and any $\ep>0$ there is $n$ such that $\fset\subset_{\ep} A_n$.
\begin{lemma}\label{lemma-1}
Let $\mathcal{C}$ be a class consisting of finite direct sums of separable
simple weakly semiprojective
 C*-algebras. Let $X$ be a compact metrizable space and let $A$ be a
   $C(X)$-algebra. Let $\fset\subset A$ be
a finite set, let $\ep>0$ and suppose  that
    $A(x)$ admits an exhaustive  sequence of
     C*-algebras isomorphic to C*-algebras in $\mathcal{C}$ for some $x \in X$.
Then there exist a compact neighborhood $U$ of $x$  and a $*$-\ho\ $\varphi: D
\to A(U)$ for some $D\in\mathcal{C}$ such that $\pi_U(\fset)\subset_\ep \varphi
(D)$. If  $A$ is a continuous $C(X)$-algebra, then we may arrange that
  $\varphi_z$ is injective for all $z\in
U$.
\end{lemma}
\begin{proof} Let $\fset=\{a_1,\dots,a_r\}$ and $\ep$ be given.
  By hypothesis
 there exist $D\in\mathcal{C}$, $\{c_1,\dots,c_r\}\subset D$ and a
 $*$-\mo\ $\iota: D \to A(x)$ such that
  $\|\pi_x(a_i)- \iota(c_i)\|<\ep/2$, for
all $i=1,...,r$.  Set $U_n=\{y\in X:d(x,y)\leq 1/n\}$. Choose a full element
$d_j$ in each direct summand of $D$.
 Since $D$ is \wsp, there
is a
 $*$-\ho\ $\varphi:D \to A(U_n)$ (for some $n$) such that $\|\pi_x\varphi
(c_{i})-\iota(c_{i})||<\ep/2$ for all $i=1,...,r$, and $\|\pi_x\varphi
(d_j)-\iota(d_j)\|\leq\|d_j\|/2$ for all $d_j$.
 Therefore
 \[\|\pi_x\varphi(c_i)-\pi_x(a_i)\|\leq \|\pi_x\varphi(c_i)-\iota(c_i)\|+
 \|\pi_x(a_i)-
\iota(c_i)\|<\ep/2+\ep/2=\ep\]  and $\varphi_x$ is injective by
Lemma~\ref{fullness}. By Lemma~\ref{a(x)delta}(i), after increasing $n$ and
setting $U=U_n$ and $\varphi=\pi_U\varphi$, we have
\[\|\varphi(c_i)-\pi_U(a_i)\|=\|\pi_U\big(\varphi(c_i)-a_i\big)\|<\ep,\]
for all $i=1,...,r$. This shows that $\pi_U(\fset)\subset_\ep \varphi (D)$.
 If $A$ is
continuous, then after shrinking $U$ we may arrange that
$\|\varphi_z(d_j)\|\geq\|\varphi_x(d_j)\|/2=\|d_j\|/2$ for all $d_j$ and all
$z\in U$. This implies that $\varphi_z$ in injective for all $z\in U$.
\end{proof}
\begin{lemma}\label{lemma-1+sp}
 Let $X$ be a compact metrizable space and let $A$ be a
   separable continuous $C(X)$-algebra the fibers of which are
   stable Kirchberg algebras. Let $\fset\subset A$ be a finite set and let $\ep>0$.
   Suppose that there exist a     KK-semiprojective stable Kirchberg algebra $D$
   and $\sigma\in KK(D,A)$
  such that $\sigma_x\in KK(D, A(x))^{-1}$
for some $x\in X$. Then there exist a closed neighborhood $U$ of $x$ and a full
$*$-\ho\ $\psi:D \to A(U)$ such that $KK(\psi)=\sigma_U$ and
$\pi_U(\fset)\subset_\ep \psi(D)$.
\end{lemma}
\begin{proof} By \cite[Thm.~8.4.1]{Ror:encyclopedia} there is an isomorphism
$\psi_0:D \to A(x)$ such that $KK(\psi_0)=\sigma_x$. Let $\mathcal{H}\subset D$
be such that  $\psi_0(\mathcal{H})=\pi_x(\fset)$. Set $U_n=\{y\in X:d(x,y)\leq
1/n\}$. By Theorem~\ref{equivalences} $D$ is KK-stable and \wsp. By
Proposition~\ref{char-semiprojectivity} there exists a $*$-\ho\ $\psi_n:D\to
A(U_n)$ (for some $n$) such that $\|\pi_x\psi_n(d)-\psi_0(d)\|<\ep$ for all
$d\in \hset$ and $KK(\pi_x\psi_n)=KK(\psi_0)=\sigma_x$. Since
$\varinjlim_{\,m}\, KK(D,A(U_m))=KK(D,A(x))$, we deduce that there is $m\geq n$
such that $KK(\pi_{U_m}\psi_n)=\sigma_{U_{m}}$. By increasing $m$ we may
arrange that $\pi_{U_m}(\fset)\subset_\ep \pi_{U_m}\psi_n(D)$ since we have
seen that $\pi_x(\fset)=\psi_0(\hset)\subset_\ep \pi_x\psi_n(D)$. We can
arrange that $\psi_z$ is injective for all $z\in U$ by reasoning as in the
proof of Lemma~\ref{lemma-1}.
 We conclude by setting $U=U_m$
and $\psi=\pi_{U_m}\psi_n$.
\end{proof}
 The following lemma is useful for  constructing  fibered morphisms.
\begin{lemma}\label{multi-perturbation-lemma}
 Let $(D_j)_{j\in J}$ be a finite family consisting of finite direct sums
 of  weakly semiprojective simple C*-algebras.
  Let $\ep>0$ and for each $j\in J$ let
$\hset_j\subset D_j$ be a finite set such that for each direct summand of $D_j$
there is an element of $\hset_j$ of norm $\geq \ep$ which is contained and is
full in that summand.
 Let $\gset_j\subset D_j$ and
$\delta_j>0$ be given by Proposition~\ref{char-semiprojectivity} applied to
$D_j$, $\hset_j$ and $\ep/2$.
    Let $X$ be
a compact metrizable space, let $(Z_j)_{j\in J}$ be  disjoint nonempty closed
subsets of $X$ and let $Y$ be a closed nonempty subset of $X$ such that
$X=Y\cup \big(\cup_j Z_j\big)$. Let $A$ be a continuous $C(X)$-algebra and let
$\fset$ be a finite subset of $A$. Let $\eta:B(Y)\rightarrow A(Y)$ be a
 $*$-monomorphism of $C(Y)$-algebras and let $\varphi_j:D_j \to A(Z_j)$ be
 $*$-\hos\ such that $(\varphi_j)_x$ is injective for all $x\in Z_j$ and
 $j\in
 J$, and which satisfy the following conditions:
\begin{itemize}
        \item[(i)] $\pi_{Z_j}(\fset)\subset_{\ep/2} \varphi_j(\hset_j)$, for all $j \in
        J$,
        \item[(ii)] $\pi_Y(\fset)\subset_{\ep} \eta(B)$,
        \item[(iii)] $\pi^{Z_j}_{Y\cap Z_j}\varphi_j(\gset_j)
        \subset_{\delta_j} \pi^Y_{Y\cap Z_j}\eta(B)$, for all $j \in J$.
    \end{itemize}
Then, there are $C(Z_j)$-linear  $*$-\mos\ $\psi_j:C(Z_j)\ot D_j\to A(Z_j)$,
satisfying \begin{equation}\label{perturb-on-H}
    \|\varphi_j(c)-\psi_j(c)\|<\ep/2, \,\,\text{for all}\,\,
    c\in \hset_j,\,\,\text{and}\,\, j\in J,
\end{equation}
and such that if we set $E=\bigoplus_j C(Z_j)\otimes D_j$, $Z=\cup_j Z_j$, and
$\psi:E \to A(Z)=\bigoplus_j A(Z_j)$, $\psi=\oplus_j\psi_j$, then $\pi^Z_{Y\cap
Z} (\psi(E))\subseteq \pi^Y_{Y\cap Z} (\eta(B))$, $\pi_Z(\fset)\subset_\ep
\psi(E)$ and hence
\[\fset \subset_{\ep}\eta(B)\oplus_{Y\cap Z} \psi(E)=\chi( B
\oplus_{\pi\eta,\pi\psi} E),
 \]
where $\chi$ is the  isomorphism   induced by the pair $(\eta,\psi)$. If we
assume that each $D_j$ is $KK$-stable, then we also have
$KK(\varphi_j)=KK(\psi_j|_{D_j})$ for all $j\in J$.\end{lemma}
\begin{proof} Let $\fset=\{a_1,\dots,a_r\}\subset A$ be as in the statement. By
 (i), for each $j\in J$ we find
 $\{c_1^{(j)},\dots,c_r^{(j)}\}\subseteq \hset_j$ such that
$\|\varphi_j(c_i^{(j)})-\pi_{Z_j}(a_i)\|<\ep/2$ for all $i$. Consider the
$C(X)$-algebra $A\oplus_Y\eta (B)\subset A$. From (iii),
 Lemma~\ref{a(x)delta}(iv) and Lemma~\ref{sumover-Y} we obtain
\begin{equation*}\label{101iii} \varphi_j(\gset_j)\subset_{\delta_j}
\pi_{Z_j}(A\oplus_Y\eta (B)).
\end{equation*}
 Applying Proposition~\ref{char-semiprojectivity} we perturb
$\varphi_j$ to a  $*$-\ho\ $\psi_j:D_j \to \pi_{Z_j}(A\oplus_Y\eta (B))$
satisfying ~\eqref{perturb-on-H}, and hence such that $\|
\varphi_j(c_i^{(j)})-\psi_j(c_i^{(j)})\|<\ep/2,$  for all $i,j$. Therefore
\begin{equation*}\label{5}
\| \psi_j(c_i^{(j)})-\pi_{Z_j}(a_i)\|\leq \|
\psi_j(c_i^{(j)})-\varphi_j(c_i^{(j)})\|+\|\varphi_j(c_i^{(j)})-\pi_{Z_j}(a_i)\|<\ep.
\end{equation*}  This
shows that $\pi_{Z_j}(\fset)\subset_{\ep}\psi_j(D_j)$. From
~\eqref{perturb-on-H}  and Lemma~\ref{fullness} we obtain that each
$(\psi_j)_x$ is injective.
 We extend $\psi_j$ to a
$C(Z_j)$-linear  $*$-\mo\ $\psi_j:C(Z_j)\otimes D_j \to
 \pi_{Z_j}(A\oplus_Y\eta (B))$ and then we define $E$, $\psi$ and $Z$
  as in the statement. In this way
we obtain that $\psi:E \to (A\oplus_Y\eta(B))(Z)\subset A(Z)$ satisfies
\begin{equation}\label{abc}
   \pi_Z(\fset)\subset_{\ep} \psi(E).
\end{equation}
The property $\psi(E)\subset (A\oplus_Y\eta(B))(Z)$ is equivalent to
$\pi^Z_{Y\cap Z} (\psi(E))\subseteq \pi^Y_{Y\cap Z} (\eta(B))$ by
Lemma~\ref{aaa}(b). Finally, from (ii), \eqref{abc} and Lemma~\ref{aaa}\,(c) we
get $\fset \subset_{\ep} \eta(B)\oplus_{Y\cap Z}
  \psi(E).$
\end{proof}
Let $\mathcal{C}$ be as in Lemma~\ref{lemma-1}. Let $A$ be a
    $C(X)$-algebra,
  let $\fset\subset A$ be a finite set and let $\ep>0$.
    An $(\fset,\ep,\mathcal{C})$-\emph{approximation} of $A$
\begin{equation}\label{FE-approx}
\alpha=\{\fset,\ep,\{U_i,\varphi_i:D_i\to
A(U_i),\mathcal{H}_i,\mathcal{G}_i,\delta_i\}_{i\in I}\},
\end{equation}
is a collection with the following properties: $(U_i)_{i\in I}$ is a finite
family of closed subsets of $X$,  whose interiors  cover  $X$ and
 $(D_i)_{i\in I}$ are C*-algebras in $\mathcal{C}$;
 for each $i\in I$,
$\varphi_i:D_i\to A(U_i)$ is a  $*$-homomorphism such that $(\varphi_i)_x$ is
injective for all $x\in U_i$;
  $\mathcal{H}_i\subset D_i$ is a finite
  set  such
that $\pi_{U_i}(\fset)\subset_{\ep/2} \varphi_i(\mathcal{H}_i)$ and such that
for each direct summand of $D_i$ there is an element of $\hset_i$ of norm $\geq
\ep$ which is contained and is full in that summand; the finite set
$\mathcal{G}_i\subset D_i$ and $\delta_i>0$ are given by
Proposition~\ref{char-semiprojectivity} applied to the weakly semiprojective
C*-algebra $D_i$ for the input data $\mathcal{H}_i$ and $\ep/2$; if  $D_i$ is
$KK$-stable, then $\mathcal{G}_i$ and $\delta_i$ are chosen such that the
second part of Proposition~\ref{char-semiprojectivity} also applies.
\begin{lemma}\label{why-approx-exists}
   Let $A$   and  $\mathcal{C}$ be as in Lemma~\ref{lemma-1}.
    Suppose that each fiber of $A$
    admits an exhaustive  sequence of
     C*-algebras isomorphic to C*-algebras in $\mathcal{C}$.
     Then for  any finite subset $\fset$ of $A$   and any $\ep>0$
     there is an $(\fset,\ep,\mathcal{C})$-approximation of $A$.
     Moreover, if $A$, $D$ and $\sigma$ are as in Lemma~\ref{lemma-1+sp}
     and $\sigma_x\in KK(D,A(x))^{-1}$ for all $x \in X$,
     then there is an $(\fset,\ep, \mathcal{C})$-approximation of $A$
       such that $\mathcal{C}=\{D\}$ and
      $KK(\varphi_i)=\sigma_{U_i}$ for
     all $i\in I$.
\end{lemma}
\begin{proof} Since $X$ is compact, this is an immediate consequence of Lemmas~\ref{lemma-1},
\ref{lemma-1+sp}
and Proposition~\ref{char-semiprojectivity}.
\end{proof}

 It is
useful to consider the following operation of restriction. Suppose that $Y$ is
a closed subspace of $X$ and let $(V_j)_{j\in J}$ be a finite family of closed
subsets of $Y$ which refines the family $(Y\cap U_i)_{i\in I}$ and such that
the interiors of the $V_j\,'s$ form a cover of $Y$. Let $\iota:J\to I$ be a map
such that $V_j\subseteq Y\cap U_{\iota(j)}$. Define
\[\iota^*(\alpha)=\{\pi_Y(\fset),\ep,\{V_j,\pi_{V_j}\varphi_{\iota(j)}:D_{\iota(j)}\to
A(V_j),\mathcal{H}_{\iota(j)},\mathcal{G}_{\iota(j)},\delta_{\iota(j)}\}_{j\in
J}\}.\] It is obvious that  $\iota^*(\alpha)$ is a
$(\pi_Y(\fset),\ep,\mathcal{C})$-approximation of $A(Y)$. The operation
$\alpha\mapsto \iota^*(\alpha)$  is useful even  in the case $X=Y$. Indeed, by
applying this procedure we can refine the cover of $X$ that appears in a given
$(\fset,\ep,\mathcal{C})$-approximation of $A$.

    An
$(\fset,\ep,\mathcal{C})$-approximation $\alpha$  (as in \eqref{FE-approx})
  is subordinated to an
$(\fset',\ep',\mathcal{C})$-approximation,
$\alpha'=\{\fset',\ep',\{U_{i'},\varphi_{i'}:D_{i'}\to
A(U_{i'}),\mathcal{H}_{i'},\mathcal{G}_{i'},\delta_{i'}\}_{{i'}\in I'}\},$
written $\alpha\prec\alpha'$, if
\begin{itemize}
        \item[(i)] $\fset\subseteq \fset'$,
        \item[(ii)] $\varphi_i(\mathcal{G}_i)\subseteq \pi_{U_i}(\fset')$
        for all $i\in I$, and
        \item[(iii)] $\ep'< \min\big(\{\ep\}\cup\{\delta_i,\, i\in I\}\big)$.
    \end{itemize}
Let us note that, with notation as above,  we have $\iota^*(\alpha)\prec
\iota^*(\alpha')$ whenever $\alpha\prec\alpha'$.

The following theorem is the crucial technical result of our paper. It provides
an approximation of continuous $C(X)$-algebras by subalgebras of category $\leq
\mathrm{dim}(X)$.
\begin{theorem}\label{basic-approx}
Let $\mathcal{C}$ be a class consisting of finite direct sums of weakly
semiprojective simple C*-algebras.
  Let $X$ be a finite dimensional compact metrizable space
    and let $A$ be a separable continuous
    $C(X)$-algebra the fibers of which admit exhaustive sequences
     of C*-algebras isomorphic to C*-algebras in $\mathcal{C}$.
 For any finite set $\fset\subset A$
 and any $\ep>0$ there exist $n\leq\mathrm{dim}(X)$ and an $n$-fibered $\mathcal{C}$-monomorphism
$(\psi_0,\dots,\psi_n)$ into $A$   which induces a $*$-\mo\
$\eta:A(\psi_0,...,\psi_n)\to A$ such that  $\fset\subset_\ep \eta
(A(\psi_0,...,\psi_n))$.
\end{theorem}
\begin{proof}
By Lemma~\ref{why-approx-exists}, for any finite set $\fset\subset A$
 and any $\ep>0$
  there is an $(\fset,\ep,\mathcal{C})$-approximation of $A$. Moreover, for
  any finite set $\fset\subset A$,
  any $\ep>0$ and any $n$, there is a sequence
  $\{\alpha_k:\, 0\leq k\leq n\}$ of $(\fset_k,\ep_k,\mathcal{C})$-approximations
  of $A$ such that $(\fset_0,\ep_0)=(\fset,\ep)$ and $\alpha_k$ is
   subordinated to $\alpha_{k+1}$:
   \[\alpha_0\prec\alpha_1\prec\cdots\prec\alpha_n.\]
Indeed, assume that $\alpha_k$ was constructed. Let us choose a finite set
$\fset_{k+1}$ which contains $\fset_k$ and liftings to $A$ of all the elements
in $\bigcup_{i_k\in I_k}\varphi_{i_k}(\gset_{i_k}).$ This choice takes care of
the above conditions (i) and (ii). Next we choose $\ep_{k+1}$ sufficiently
small such that (iii) is satisfied. Let $\alpha_{k+1}$ be an
$(\fset_{k+1},\ep_{k+1},\mathcal{C})$-approximation  of $A$ given by
Lemma~\ref{why-approx-exists}. Then obviously $\alpha_k\prec\alpha_{k+1}$. Fix
a tower of approximations of $A$ as above where $n=\mathrm{dim}(X)$.

By \cite[Lemma~3.2]{BK:bundles}, for every open cover $\mathcal{V}$ of $X$
there is a finite open cover $\mathcal{U}$ which refines $\mathcal{V}$ and such
that the set $\mathcal{U}$ can be partitioned into $n+1$ nonempty subsets
consisting of elements with pairwise disjoint closures.
 Since we can refine simultaneously
the covers that appear in a finite family $\{\alpha_k:\, 0\leq k\leq n\}$ of
approximations while preserving subordination,
 we may arrange not only that all $\alpha_k$ share the same cover $(U_i)_{\in I},$
but moreover, that the cover $(U_i)_{i\in I}$ can be partitioned into $n+1$
subsets $\U_0,\dots,\U_n$ consisting of mutually disjoint elements. For
definiteness, let us write $\U_k=\{U_{i_k}:i_k\in I_k\}$. Now for each $k$ we
consider the closed subset of $X$
\[Y_k=\bigcup_{i_k\in I_k}\,U_{i_k},\]
the map $\iota_k:I_k\to I$ and the
$(\pi_{Y_{k}}(\fset_k),\ep_k,\mathcal{C})$-approximation of $A(Y_k)$, induced
by $\alpha_k$,  which is of the form
\[ \iota_k^*(\alpha_k)=\{\pi_{Y_k}(\fset_k),\ep,\{U_{i_k},\varphi_{i_k}:D_{i_k}\to
A(U_{i_k}),\mathcal{H}_{i_k},\mathcal{G}_{i_k},\delta_{i_k}\}_{{i_k}\in
I_k}\},\] where each $U_{i_k}$ is nonempty.
 We have
\begin{equation}\label{0}
   \pi_{U_{i_k}}(\fset_k)\subset_{\ep_k/2}\varphi_{i_k}(\mathcal{H}_{i_k}),
\end{equation}
 by construction. Since $\alpha_k\prec\alpha_{k+1}$ we obtain
\begin{equation}\label{1}
    \fset_k\subseteq \fset_{k+1},
\end{equation}
\begin{equation}\label{2}
    \varphi_{i_k}(\gset_{i_k})\subseteq\pi_{U_{i_k}}(\fset_{k+1}),\,\,\text{for
    all}\,\, i_k\in I_k,
\end{equation}
\begin{equation}\label{3}
   \ep_{k+1}< \min\big(\{\ep_k\}\cup\{\delta_{i_k},\, {i_k}\in I_k\}.\big)
\end{equation}
Set $X_k=Y_k\cup \dots \cup Y_n$ and $E_k=\oplus_{i_k} \,C(U_{i_k})\ot D_{i_k}$
for $0\leq k \leq n$.
 We shall construct  a sequence  of $C(Y_k)$-linear $*$-\mos, $\psi_k:E_k \to
A(Y_k)$, $k=n,...,0$,  such that $(\psi_k,\dots,\psi_n)$ is an $(n-k)$-fibered
monomorphism into $A(X_k)$. Each map
\[\psi_k=\oplus_{i_k}\psi_{i_k}:E_k \to A(Y_k)=\oplus_{i_k}\,A(U_{i_k})\]
will have components
 $\psi_{i_k}:C(U_{i_k})\ot D_{i_k}\to
A(U_{i_k})$  whose restrictions to $D_{i_k}$ will be perturbations of
$\varphi_{i_k}:D_{i_k}\to A(U_{i_k})$, ${i_k}\in I_k$.
 We shall construct the maps $\psi_k$  by induction on decreasing $k$
  such that if
 $ B_{k}=A(X_k)(\psi_k,\dots,\psi_n)$ and  $\eta_{k}:B_{k}\to A(X_k)$
 is the map induced by the
 $(n-k)$-fibered monomorphism $(\psi_k,\dots,\psi_n)$, then
\begin{equation}\label{6k}
    \pi_{X_{k+1}\cap U_{i_k}}\big(\psi_{i_k}(D_{i_k})\big)
    \subset \pi_{X_{k+1}\cap U_{i_k}}
    \big(\eta_{k+1}(B_{k+1})\big), \,\forall \,{i_k} \in I_k,
\end{equation} and
\begin{equation}\label{7k}
    \pi_{X_k}(\fset_k)\subset_{\ep_k}\eta_{k}(B_{k}).
\end{equation}
Note that ~\eqref{6k} is equivalent to
\begin{equation}\label{6kk}
    \pi_{X_{k+1}\cap Y_k}\big(\psi_k(E_k)\big)\subset \pi_{X_{k+1}\cap Y_k}
    \big(\eta_{k+1}(B_{k+1})\big).
\end{equation}
 For the first step of
induction, $k=n$, we choose  $\psi_n=\oplus_{i_n}\widetilde{\varphi}_{i_n}$
where $\widetilde{\varphi}_{i_n}:C(U_{i_n})\ot D_{i_n}\to A(U_{i_n})$ are
$C(U_{i_n})$-linear extensions of the original $\varphi_{i_n}$. Then $B_n= E_n$
and $\eta_n=\psi_n$.
 Assume  that
$\psi_n,\dots,\psi_{k+1}$ were constructed and that they have the desired
properties. We shall construct now $\psi_{k}$. Condition~\eqref{7k} formulated
for $k+1$ becomes
\begin{equation}\label{7kk}
    \pi_{X_{k+1}}(\fset_{k+1})\subset_{\ep_{k+1}}\eta_{k+1}(B_{k+1}).
\end{equation}
  Since $\ep_{k+1}<\delta_{i_k}$, by using  \eqref{2} and
~\eqref{7kk}  we obtain
\begin{equation}\label{iii}
\pi_{X_{k+1}\cap
U_{i_k}}\big(\varphi_{i_k}(\gset_{i_k})\big)\subset_{\delta_{i_k}}
\pi_{X_{k+1}\cap U_{i_k}}
    \big(\eta_{k+1}(B_{k+1})\big), \,\,\text{for all}\, \,{i_k} \in I_k.
\end{equation}
  Since $\fset_k\subseteq \fset_{k+1}$
 and $\ep_{k+1}<{\ep_k}$, condition~\eqref{7kk} gives
\begin{equation}\label{ii}
    \pi_{X_{k+1}}(\fset_{k})\subset_{\ep_k}\eta_{k+1}(B_{k+1}).
\end{equation}
 Conditions \eqref{0},
\eqref{iii}  and  \eqref{ii} enable us to apply
Lemma~\ref{multi-perturbation-lemma} and perturb $\widetilde\varphi_{i_k}$ to a
$*$-\mo\ $\psi_{i_k}:C(U_{i_k})\ot D_{i_k}\to A(U_{i_k})$ satisfying
 ~\eqref{6k} and ~\eqref{7k} and such that
 \begin{equation}\label{perturb-with-homotopy-class}
    KK(\psi_{i_k}|_{D_{i_k}})=KK(\varphi_{i_k})
\end{equation}
if the algebras in $\mathcal{C}$ are assumed to be $KK$-stable. We set
 $\psi_k=\oplus_{i_k}\,\psi_{i_k}$ and
  this completes the construction of $(\psi_0,\dots,\psi_n)$.
  Condition ~\eqref{7k} for $k=0$ gives
$\fset\subset_\ep \eta_0(B_0)=\eta(A(\psi_0,\dots,\psi_n)).$ Thus
$(\psi_0,\dots,\psi_n)$ satisfies the conclusion of the theorem. Finally let us
note that it can happen that $X_k=X$ for some $k>0$. In this case
$\fset\subset_\ep A(\psi_k,...,\psi_n)$ and for this reason we write $n\leq
\mathrm{dim}(X)$ in the statement of the theorem.
\end{proof}

\begin{proposition}\label{approximation-fiber-is-D}
    Let $X$ be a finite dimensional
 compact metrizable space and
 let $A$ be a separable
continuous $C(X)$-algebra the fibers of which are
 stable Kirchberg algebras.
Let $D$ be a  $KK$-semiprojective stable Kirchberg
   algebra and suppose that  there exists $\sigma\in KK(D,A)$
   such that  $\sigma_x\in
 KK(D,A(x))^{-1}$ for all $x\in X$.
For any finite subset $\fset$ of $A$
 and any $\ep>0$ there is an $n$-fibered $\mathcal{C}$-monomorphism
 $(\psi_0,\dots,\psi_n)$ into $A$
  such that $n\leq\mathrm{dim}(X)$, $\mathcal{C}=\{D\}$,
 and each component $\psi_i:C(Y_i)\ot D \to A(Y_i)$ satisfies
 $KK(\psi_i)=\sigma_{Y_i}$,
 $i=0,\dots,n.$ Moreover, if $\eta:A(\psi_0,\dots,\psi_n)\to A$ is
 the  induced $*$-\mo,
 then $\fset\subset_\ep \eta(A(\psi_0,\dots,\psi_n))$ and $KK(\eta_x)$ is a
 KK-equivalence for each $x\in X$.
\end{proposition}
\begin{proof}
 We repeat the proof of
Theorem~\ref{basic-approx} while using only
$(\fset_i,\ep_i,\{D\})$-approximations of $A$ provided by the second part of
Lemma~\ref{why-approx-exists}. The outcome will be an $n$-fibered
$\{D\}$-monomorphism $(\psi_0,\dots,\psi_n)$ into $A$
  such
 that $\fset\subset_\ep \eta(A(\psi_0,\dots,\psi_n))$.
 Moreover we can arrange that
 $KK(\psi_i)=\sigma_{Y_i}$ for all $i=0,\dots,n$,
 by~\eqref{perturb-with-homotopy-class},
  since $KK(\varphi_{i_k})=\sigma_{U_{i_k}}$ by Lemma~\ref{why-approx-exists}.
   If $x\in X$, and $i=\min\{k:x\in Y_k\}$, then
  $\eta_x\equiv(\psi_i)_x$, and hence $KK(\eta_x)$ is a KK-equivalence.
\end{proof}
\begin{remark}\label{polys} Let us point out that we can strengthen the conclusion
of Theorem~\ref{basic-approx} and Proposition~\ref{approximation-fiber-is-D} as
follows. Fix a metric $d$ for the topology of $X$. Then we may arrange that
there is a closed cover $\{Y'_0,...,Y'_n\}$ of $X$
 and a number $\ell>0$ such that
 $\{x:d(x,Y'_i)\leq \ell\}\subset Y_i$ for $i=0,...,n$.
 Indeed, when we choose the finite closed cover $\mathcal{U}=(U_i)_{i\in I}$ of $X$
 in the proof of
 Theorem~\ref{basic-approx} which can be partitioned into $n+1$
subsets $\U_0,\dots,\U_n$ consisting of mutually disjoints elements,
  as given by \cite[Lemma~3.2]{BK:bundles},
 and which refines all the covers
  $\mathcal{U}(\alpha_0),...,\mathcal{U}(\alpha_n)$
  corresponding to $\alpha_0,...,\alpha_n$,
 we may assume that $\mathcal{U}$ also refines the covers given by the
 interiors of  the elements of
 $\mathcal{U}(\alpha_0),...,\mathcal{U}(\alpha_n)$.
 Since each $U_i$ is compact and $I$ is finite, there is $\ell>0$ such that
 if $V_i=\{x:d(x,U_i)\leq \ell\}$, then the cover $\mathcal{V}=(V_i)_{i\in I}$
 still refines all of $\mathcal{U}(\alpha_0),...,\mathcal{U}(\alpha_n)$
 and for each $k=0,...,n$, the elements
 of $\mathcal{V}_k=\{V_i:U_i\in \mathcal{U}_k\}$,
 are still mutually
 disjoint. We shall use the cover
 $\mathcal{V}$ rather than $\mathcal{U}$ in the proof of the two theorems
 and observe that $Y'_k\stackrel{def}{=}\bigcup_{i_k\in I_k} U_{i_k}\subset
 \bigcup_{i_k\in I_k} V_{i_k}=Y_k$ has the desired property.
 Finally let us note that if we define $\psi_i':E(Y_i')\to A(Y_i')$
 by $\psi_i'=\pi_{Y_i'}\psi_i$, then $(\psi_0',\dots,\psi_n')$
 is an $n$-fibered $\mathcal{C}$-monomorphism into $A$ which satisfies
 the conclusion of Theorem~\ref{basic-approx} and
  Proposition~\ref{approximation-fiber-is-D} since $\pi_{Y'_i}(\fset)\subset_\ep \psi_i'(E_i)$
  for all $i=0,\dots,n$ and $X=\bigcup_{i=1}^{\,n}Y_i'$.
\end{remark}
\section{Representing $C(X)$-algebras as inductive limits}\label{three}
 We have seen that Theorem~\ref{basic-approx}
 yields exhaustive sequences for certain  $C(X)$-algebras. In this
section we show how to pass from an exhaustive sequence to a nested exhaustive
sequence using semiprojectivity. The remainder of the paper does not  depend on
this section.
\begin{proposition}\label{basic-approx-proj}
Let $X$, $A$ and  $\mathcal{C}$ be as in Theorem~\ref{basic-approx}. Let
$(\psi_0,\dots,\psi_n)$ be an $n$-fibered $\mathcal{C}$-monomorphism into $A$
with components $\psi_i:E_i \to A(Y_i)$.
 Let $\fset_i\subset E_i$, $\fset\subset A(\psi_0,...,\psi_n)$ be finite sets
and let $\ep>0$. Then there are finite sets $\gset_i\subset E_i$ and
$\delta_i>0$, $i=0,...,n$,  such that for any $C(X)$-subalgebra $A' \subset A$
which satisfies $\psi_i(\gset_i)\subset_{\delta_i} A'(Y_i)$, $i=0,...,n$, there
is an $n$-fibered $\mathcal{C}$-monomorphism $(\psi_0',\dots,\psi_n')$ into
$A'$, with $\psi'_i:E_i \to A'(Y_i)$ and such that (i)
$\|\psi_i(a)-\psi'_i(a)\|<\ep$ for all $a\in \fset_i$ and all $i\in
\{0,...,n\}$, (ii) $(\psi_j)_x^{-1}(\psi_i)_x=(\psi'_j)_x^{-1}(\psi'_i)_x$ for
all $x\in Y_i\cap Y_j$ and $0\leq i\leq j\leq n$. Moreover
$A(\psi_0,\dots,\psi_n)= A'(\psi_0',\dots,\psi_n')$ and
 the maps
 $\eta:A(\psi_0,...,\psi_n)\to A$ and $\eta':A'(\psi_0',\dots,\psi_n')\to A'$
 induced by $(\psi_0,\dots,\psi_n)$ and $(\psi_0',\dots,\psi_n')$ satisfy
(iii) $\|\eta(a)-\eta'(a)\|<\ep$ for all $a \in \fset$.
\end{proposition}
\begin{proof}  Let us observe that if we prove (i) and (ii) 
then (iii) will follow by enlarging the sets $\fset_i$ so that $p_i(\fset)\subset \fset_i$,
where $p_i:A(\psi_0,...,\psi_n)\to E_i$
are the coordinate maps.
We proceed now with the proof of (i) and (ii) by making some simplifications.
 We may assume that $E_0=C(Y_0)\ot D_0$ with $D_0\in
\mathcal{C}$ since the perturbations corresponding to disjoint closed sets can
be done independently of each other. Without any loss of generality, we may assume that
$\fset_0\subset D_0$ since we are working with morphisms on $E_0$ which are
$C(Y_0)$-linear.
We also enlarge $\fset_0$
 so that for each direct summand $C$ of $D_0$, $\fset_0$ contains an
  element $c$ which is full in  $C$ and such that $\|c\|\geq 2\ep$.

 The proof is by induction on $n$.
 If
$n=0$ the statement follows from Proposition~\ref{char-semiprojectivity} and
Lemma~\ref{fullness}. Assume now that the statement is true for $n-1$. Let
$E_i$,  $\psi_i$, $A$, $A'$, $\fset_i$, $1\leq i \leq n$  and $\ep$ be as in the
statement. For $0\leq i<j \leq n$ let $\eta_{j,i} :E_i(Y_i\cap Y_j)\to E_j(Y_i\cap Y_j)$ be the $*$-\ho\
of $C(Y_i\cap Y_j)$-algebras defined fiberwise by $(\eta_{j,i})_x=(\psi_j)^{-1}_{x}(\psi_i)_{x}$

Let $\gset_0$ and $\delta_0$ be given by
Proposition~\ref{semiproj_combined} applied to the C*-algebra $D_0$ for the
input data $\fset_0$ and $\ep$. For each $1\leq j \leq n$ choose a finite subset $\hset_j$
of $E_j$ whose restriction to $Y_j\cap Y_0$ contains $\eta_{j,0}(\gset_0)$.
Consider the sets
$\fset_j^{\prime}:=\fset_j\cup \hset_j$, $1 \leq j \leq n$ and the number $\ep'=\min\{\delta_0,\ep\}$.
Let $\gset_1,...\gset_n$ and
$\delta_1,...,\delta_n$ be given by the inductive assumption for $n-1$ applied
to $A(X_1)$, $A'(X_1)$, $\psi_j$, $\fset^{\prime}_j$, $1 \leq j \leq n$ and $\ep'$, where $X_1=Y_1\cup\cdots\cup Y_n$.

  We need to show that $\gset_0,\gset_1,...\gset_n$ and
  $\delta_0,\delta_1,...,\delta_n$ satisfy the statement.
 By
  the inductive step there exists an $(n-1)$-fibered $\mathcal{C}$-monomorphism
  $(\psi'_1,\dots,\psi'_n)$ into $A'(X_1)$ with components $\psi_j':E_j\to A'(Y_j)$
    such that

(a)\quad $\|\psi_j(a)-\psi'_j(a)\|<\ep'=\min\{\delta_0,\ep\}$ for all $a\in \fset_j\cup \hset_j$ and all
$1\leq j \leq n$,

(b)\quad $(\psi_j)_x^{-1}(\psi_i)_x=(\psi'_j)_x^{-1}(\psi'_i)_x$ for all $x\in
Y_i\cap Y_j$ and $1\leq i\leq j\leq n$,

The condition (b) enables to define a $*$-\ho\ $\varphi:E_0 \to A'(Y_0\cap X_1 )$
with fiber maps  $\varphi_x=(\psi_j')_x(\psi_j)^{-1}_x (\psi_0)_x $ for $x\in Y_0\cap Y_j $ and $1\leq j\leq n$.

Let us observe that $\psi_0:E_0 \to A(Y_0)$ is an approximate lifting of $\varphi$.
More precisely we have $\|\pi_{X_1\cap Y_0}^{Y_0}\psi_0(a)-\varphi(a)\|<\delta_0$ for all
$a\in \gset_0$. Indeed, for $x\in Y_0\cap Y_j$, $1\leq j \leq n$ and $a\in \gset_0$ we have
\begin{eqnarray*}
 \|(\psi_0)_x (a(x))-(\psi'_j)_x(\psi^{-1}_j)_x(\psi_0)_x(a(x))\|&=&
\|(\psi_j)_x(\eta_{j,0})_x(a(x))-(\psi'_j)_x(\eta_{j,0})_x(a(x))\| \\&\leq&
\sup_{h\in \hset_j}\|\psi_j(h)-\psi'_j(h)\|<\ep'\leq \delta_0.
\end{eqnarray*}
Since we also have $\psi_0(\gset_0)\subset_{\delta_0} A'(Y_0)$ by hypothesis, it follows from Proposition~\ref{semiproj_combined} that there exists $\psi'_0:D_0 \to A(Y_0)$ such that
$\|\psi'_0(a)-\psi_0(a)\|<\ep$ for all $a\in \fset_0$ and $\pi_{Y_0\cap X_1}^{Y_0}\psi'_0=\varphi$.  By Lemma~\ref{fullness} each $(\psi'_0)_x$ is
injective since each $(\psi_0)_x$ is injective. The $C(Y_0)$-linear
extension of $\psi'_0$ to $E_0$ satisfies $(\psi_j)_x^{-1}(\psi_0)_x=(\psi'_j)_x^{-1}(\psi'_0)_x$ for all $x\in
Y_0\cap Y_j$ and $1\leq  j\leq n$ and this completes the proof of (ii). Condition (i) follows from (b).
\end{proof}
The following result gives an inductive limit representation for
continuous $C(X)$-algebras whose fibers are inductive limits of finite direct sums of
simple semiprojective C*-algebras.  For example the fibers can be arbitrary
 AF-algebras or  Kirchberg
algebras which satisfy the UCT and whose $K_1$-groups are torsion free. Indeed,
by \cite[Prop.~8.4.13]{Ror:encyclopedia},
 these algebras are isomorphic to  inductive
limits  of sequences of  Kirchberg  algebras $(D_n)$ with finitely generated
K-theory groups and torsion free $K_1$-groups. The algebras $D_n$ are
semiprojective by \cite{Spi:semiproj}.
\begin{theorem}\label{basic-approx-K1=torsion-free}
Let $\mathcal{C}$ be a class consisting of finite direct sums of semiprojective
simple C*-algebras.
  Let $X$ be a finite dimensional compact metrizable space
  and let $A$ be a separable
   continuous $C(X)$-algebra such that all its fibers
    admit exhaustive sequences
    consisting of C*-algebras  isomorphic to
     C*-algebras in $\mathcal{C}$.
    Then $A$ is isomorphic to the inductive
    limit of a sequence of continuous $C(X)$-algebras
    $A_k$  such that
$\mathrm{cat}_\mathcal{C}(A_k) \leq \mathrm{dim}(X)$.
\end{theorem}
\begin{proof}
   By Theorem
~\ref{basic-approx} and Proposition~\ref{basic-approx-proj} we find a sequence
$(\psi_0^{(k)},...,\psi_n^{(k)})$ of $n$-fibered $\mathcal{C}$-monomorphisms
into $A$ which induces $*$-\mos\
$\eta^{(k)}:A_k=A(\psi_0^{(k)},...,\psi_n^{(k)})\to A$ with the following
properties. There is a sequence of finite sets $\fset_k\subset A_k$ and a
sequence of $C(X)$-linear $*$-\mos\ $\mu_k:A_k\to A_{k+1}$ such that
\begin{itemize}
        \item[(i)] $\|\eta^{(k+1)}\mu_k(a)-\eta^{(k)}(a)\|<2^{-k}$ for
        all $a \in \fset_k$ and all $k\geq 1$,
        \item[(ii)] $\mu_k(\fset_k)\subset \fset_{k+1}$ for all $k\geq 1$,
        \item[(iii)] $\bigcup_{j=k+1}^\infty\,
        \big(\mu_{j-1}\circ\cdots\circ\mu_{k}\big)^{-1}(\fset_j)$ is dense in $A_k$
        and
        $\bigcup_{j=k}^\infty\,\eta^{(j)}(\fset_j)$ is dense in $A$ for all $k\geq 1$.
    \end{itemize}
    Arguing as in the proof of \cite[Prop.~2.3.2]{Ror:encyclopedia},
one verifies that $$\varphi_k(a)=\lim_{j \to
\infty}\eta^{(j)}\circ\big(\mu_{j-1}\circ\cdots\circ\mu_{k}\big) (a)$$ defines
a sequence of $*$-\mos\  $\varphi_k:A_k \to A$  such that
$\varphi_{k+1}\mu_k=\varphi_k$ and the induced map   $\varphi:\varinjlim_k
(A_k,\mu_k)\to A$ is an isomorphism of $C(X)$-algebras.
\end{proof}
\begin{remark}\label{units}
    By similar arguments one proves a unital version of
    Theorem~\ref{basic-approx-K1=torsion-free}.
\end{remark}
\section{When is a fibered product locally trivial}\label{LTCFP}

 For  C*-algebras $A$, $B$ we endow the space
$\mathrm{Hom}(A,B)$ of  $*$-\hos\ with the point-norm topology. If $X$ is a
compact Hausdorff space, then $\mathrm{Hom}(A,C(X)\ot B)$ is homeomorphic to
the space of continuous maps from $X$ to $\mathrm{Hom}(A,B)$ endowed with the
compact-open topology. We shall identify a $*$-\ho\
$\varphi\in\mathrm{Hom}(A,C(X)\ot B)$ with the corresponding continuous map
$X\to \mathrm{Hom}(A,B)$, $x\mapsto \varphi_x$, $\varphi_x(a)=\varphi(a)(x)$
for all $x \in X$ and $a\in A$.
  Let $D$ be a C*-algebra and let $A$ be
a $C(X)$-algebra.  If $\alpha:D \to A$ is a $*$-\ho, let us denote by
$\widetilde{\alpha}:C(X)\ot D\to A$ its (unique) $C(X)$-linear
 extension
 and write $\widetilde{\alpha}\in \mathrm{Hom}_{C(X)}(C(X)\ot D,A)$.
 For C*-algebras $D$, $B$ we shall make without further
 comment the following identifications
 \[\mathrm{Hom}_{C(X)}(C(X)\ot D, C(X)\ot B)\equiv
 \mathrm{Hom}( D, C(X)\ot B)\equiv C(X, \mathrm{Hom}( D,  B)).\]
For a  C*-algebra $D$
 we denote by $\End(D)$ the set of full (\emph{and unital if $D$ is unital}) $*$-endomorphisms
 of $D$ and by $\End(D)^0$ the
path component of $\id_D$ in $\End(D)$. Let us consider
\begin{equation*}\label{end*}
\End(D)^*=\{\gamma\in \End(D):KK(\gamma)\in KK(D,D)^{-1}\}.
\end{equation*}
\begin{proposition}\label{Pb}
    Let  $X$ be a compact metrizable space and let $D$ be a
$KK$-semiprojective Kirchberg algebra.
    Let $\alpha:D\to C(X)\ot D$ be
    a full  (and unital, if $D$ is unital) $*$-\ho\ such that $KK(\alpha_x)\in KK(D,D)^{-1}$ for all $x\in X$.
    Then there is a full  $*$-\ho\ $\Phi:D\to C(X\times [0,1])\ot D$
    such that $\Phi_{(x,0)}=\alpha_x$ and
     $\Phi_{(x,t)}\in \Aut(D)$ for all $x\in X$ and $t\in (0,1]$.
     Moreover, if
     $\Phi_{1}:D\to C(X)\ot D$ is defined by $\Phi_{1}(d)(x)=\Phi_{(x,1)}(d)$,
     for all $ d\in D$ and $x\in X$, then $\alpha\approx_{uh} \Phi_{1}$.
\end{proposition}
\begin{proof} Since $X$ is a metrizable compact space,
 $X$ is homeomorphic to the projective limit of a sequence
of finite simplicial complexes $(X_i)$ by \cite[Thm.~10.1,
p.284]{Eil-Steen:at}. Since $D$ is $KK$-semiprojective,
$KK(D,\varinjlim\,C(X_i)\ot D)=KK(D,C(X)\ot D)$ by Theorem~\ref{equivalences}.
By Theorem~\ref{Phillips-basic}, there is $i$ and a full (and unital if $D$ is
unital) $*$-\ho\ $\varphi:D\to C(X_i)\ot D$ whose KK-class maps to
$KK(\alpha)\in KK(D,C(X)\ot D)$. To summarize, we have found a finite
simplicial complex $Y$, a continuous map $h:X\to Y$ and a continuous map
$y\mapsto \varphi_y\in \End(D)$, defined on $Y$, such that the full (and unital
if $D$ is unital) $*$-\ho\ $h^*\varphi:D\to C(X)\ot D$ corresponding to the
continuous map $x\mapsto \varphi_{h(x)}$ satisfies $KK(h^*\varphi)=KK(\alpha)$.
We may arrange that $h(X)$ intersects all the path components of $Y$ by
dropping  the  path components which are not intersected. Since $\alpha_x\in
\End(D)^*$ by hypothesis, and since $KK(\alpha_x)=KK(\varphi_{h(x)}),$ we infer
that $\varphi_y\in \End(D)^*$ for all $y\in Y$. We shall find a continuous map
$y\mapsto\psi_y\in\End(D)^*$ defined on $Y$, such that the maps $y\mapsto
\psi_y\varphi_y$ and $y\mapsto \varphi_y\psi_y$ are homotopic to the constant
map $\iota$ that takes $Y$ to $\id_D$. It is clear that it suffices to deal
separately with each path component of $Y$, so that for this part of the proof
 we may assume that $Y$
is connected. Fix a point $z\in Y$.
 By \cite[Thm.~8.4.1]{Ror:encyclopedia} there
is
 $\nu\in \Aut(D)$ such that $KK(\nu^{-1})=KK(\varphi_z)$ and
  hence $KK(\nu\varphi_z)=KK(\id_D)$. By
Theorem~\ref{Phillips-basic}, there is a unitary $u\in M(D)$ such that
$u\nu\varphi_z(-) u^*$ is homotopic to $\id_D$. Let us set
$\theta=u\nu(-)u^*\in \Aut(D)$ and observe that $\theta\varphi_z \in
\End(D)^{0}$. Since $Y$ is path connected, it follows  that the entire image of
the map $y\mapsto
 \theta\varphi_y $ is contained in $\End(D)^0$.
 Since $\End(D)^0$ is a path connected H-space with unit element,
it follows by \cite[Thm.~2.4, p462]{Whitehead:top-book} that the homotopy
classes $[Y,\End(D)^0]$ (with no condition on basepoints, since the action of
the fundamental group $\pi_1(\End(D)^0,\id_D)$ is trivial by \cite[3.6,
p166]{Whitehead:top-book}) form a group under the natural multiplication.
Therefore we find $y\mapsto \psi'_y\in \End(D)^0$ such that $y\mapsto
\psi'_y\theta\varphi_y$ and $y\mapsto \theta\varphi_y\psi'_y$ are homotopic to
$\iota$. It follows that $y\mapsto\psi_y\stackrel{def}=\psi_y'\theta$ is the
homotopic inverse of
 $y\mapsto\varphi_y$ in $[Y,\End(D)^*]$.
  Composing with $h$ we obtain that the maps $x\mapsto
\varphi_{h(x)}\psi_{h(x)}$ and $x\mapsto \psi_{h(x)}\varphi_{h(x)}$ are
homotopic to
 the constant map
that takes $X$ to $\id_D$.
  By the homotopy  invariance of
$KK$-theory we obtain that
$$KK(\widetilde {h^*\varphi} \, h^*\psi)=
KK(\widetilde {h^*\psi}  \,h^*\varphi)=KK(\iota_D),$$ where  $\widetilde
{h^*\varphi}$ and $\widetilde {h^*\psi}$ denote the $C(X)$-linear extensions of
the corresponding maps  and $\iota_D:D \to C(X)\ot D$ is defined by
$\iota_D(d)=1_{C(X)}\ot d$ for all $d\in D$.    Let us recall that
$KK(h^*\varphi)=KK(\alpha)$ and hence $KK(\widetilde
{h^*\varphi})=KK(\widetilde\alpha)$. If we set $\Psi={h^*\psi}$, then
$$KK(\widetilde {\alpha}\,  \Psi)=
KK(\widetilde {\Psi} \, \alpha)=KK(\iota_D).
$$ By
Theorem~\ref{Phillips-basic} $\widetilde {\alpha} \,\Psi\approx_u \iota_D$ and
$\widetilde {\Psi} \,\alpha\approx_u \iota_D,$ and hence
 $\widetilde {\alpha} \, \widetilde\Psi\approx_u
\mathrm{id}_{C(X)\ot D}$ and $\widetilde {\Psi}\, \widetilde\alpha\approx_u
\mathrm{id}_{C(X)\ot D}.$ By \cite[Cor.~2.3.4]{Ror:encyclopedia}, there is an
isomorphism $\Gamma:C(X)\ot D \to C(X)\ot D$ such that
$\Gamma\approx_u\widetilde\alpha$. In particular $\Gamma$ is $C(X)$-linear and
$\Gamma_x\in \Aut(D)$ for all $x\in X$.
 Replacing
$\Gamma$ by $u\Gamma( \cdot)u^*$ for some unitary $u \in M(C(X)\ot D)$ we can
arrange that $\Gamma|_D$  is arbitrarily close to $\alpha$. Therefore
$KK(\Gamma|_D)=KK(\alpha)$
 since $D$ is KK-stable. By Theorem~\ref{Phillips-basic} there is
 a
continuous map $(0,1]\to U(M(C(X)\ot D))$, $t\mapsto u_t$, with the property
that
\begin{equation*}\label{d(Y)}
\lim_{t \to 0} \|u_t\Gamma(a)u_t^*-\alpha(a)\|=0,\,\,\text{for all}\quad a \in
 D.
\end{equation*}
Therefore the equation \[\Phi_{(x,t)}=
 \left\{
\begin{array}{ll}
        \alpha_x, & \hbox{if $t=0$,} \\
        u_t(x)\Gamma_xu_t(x)^*, &\hbox{if $t\in (0,1]$,} \\
\end{array}
\right.\] defines a continuous map $\Phi:X\times [0,1]\to \End(D)^*$ which
extends $\alpha$ and such that $\Phi(X\times (0,1])\subset \Aut(D)$. Since
$\alpha$ is homotopic to $\Phi_1$, we have that $\alpha\approx_{uh} \Phi_{1}$
by Theorem~\ref{Phillips-basic}.
\end{proof}

\begin{proposition}\label{101}
    Let  $X$ be a compact metrizable space and let $D$ be a
$KK$-semiprojective Kirchberg algebra. Let $Y$ be a closed subset of
    $X$.
    Assume that a  map $\gamma:Y \to \End(D)^*$
   extends to a continuous map $\alpha:X \to \End(D)^*$.
     Then there is a continuous extension $\eta:X \to \End(D)^*$ of $\gamma$,
    such that
    $\eta(X\setminus Y)\subset \Aut(D)$.
\end{proposition}
\begin{proof}
Since the map $x\to\alpha_x$ takes values in $\End(D)^*$, by
Proposition~\ref{Pb} there exists a continuous map $\Phi:X\times [0,1]\to
\End(D)^*$ which extends $\alpha$ and such that $\Phi(X\times (0,1])\subset
\Aut(D)$. Let $d$ be a metric for the topology of $X$ such that
$\mathrm{diam}(X)\leq 1$. The equation $\eta(x)=\Phi(x,d(x,Y))$ defines a map
on $X$ that satisfies the conclusion of the proposition.
\end{proof}
\begin{lemma}\label{102}
  Let  $X$ be a compact metrizable space and let $D$ be a
$KK$-semiprojective Kirchberg algebra. Let $Y$ be a closed subset of
    $X$.
     Let $\alpha:Y\times [0,1]\cup X \times \{0\} \to \End(D)$
    be a continuous map such that
    $\alpha_{(x,0)}\in \End(D)^*$ for all $x \in X$.
Suppose that there is an open set $V$ in $X$ which contains $Y$ and such that
$\alpha$ extends to a continuous map $\alpha_V: V\times [0,1]\cup X \times
\{0\}\to \End(D)$.
    Then there is
    $\eta:X\times [0,1] \to \End(D)^*$ such that $\eta$ extends $\alpha$
    and $\eta_{(x,t)}\in \Aut(D)$ for all $x\in X\setminus Y$ and $t\in (0,1]$.
\end{lemma}

\begin{proof}
By Proposition~\ref{101} it suffices to find a continuous map
$\widehat\alpha:X\times [0,1]\to \End(D)^*$ which extends $\alpha$.   Fix a
metric $d$ for the topology of $X$ and define $\lambda:X \to [0,1]$ by
$\lambda(x)=d(x,X\setminus V)\big(d(x,X\setminus V)+d(x,Y)\big)^{-1}$. Let us
define $\widehat\alpha:X\times [0,1]\to \End(D)$ by
$\widehat\alpha_{(x,t)}=\alpha_V(x,\lambda(x)t)$ and observe that
$\widehat\alpha$ extends $\alpha$. Finally, since $\widehat\alpha_{(x,t)}$ is
homotopic to $\widehat\alpha_{(x,0)}=\alpha_{(x,0)}$, we conclude that the
image of $\widehat\alpha$ in contained in $\End(D)^*$.
\end{proof}
\begin{proposition}\label{loc-triv+section=triv}
    Let  $X$ be a compact metrizable space and let $D$ be a
$KK$-semiprojective stable  Kirchberg algebra. Let $A$ be a separable
$C(X)$-algebra which is locally isomorphic to $C(X)\ot D$. Suppose that there
is $\sigma\in KK(D,A)$ such that $\sigma_x\in KK(D, A(x))^{-1}$ for all $x \in
X$. Then there is an isomorphism of $C(X)$-algebras $\psi:C(X)\ot D \to A$ such
that $KK(\psi|_D)=\sigma$.
\end{proposition}
\begin{proof} Since $X$ is compact and $A$ is locally trivial it follows that
$\mathrm{cat}_{\{D\}}(A)<\infty$. By Lemma~\ref{lemma-needed-properties-for-B},
 $A \cong pAp\ot \OOO_\infty \ot \mathcal{K}$ for some projection $p\in A$.
 By Theorem~\ref{Phillips-basic}, there is a full $*$-\ho\ $\varphi:D\to A$
 such that $KK(\varphi)=\sigma$. We shall construct an isomorphism of $C(X)$-algebras $\psi:C(X)\ot D \to A$
such that $\psi$ is homotopic to $\widetilde\varphi$, the $C(X)$-linear
extension of $\varphi$. Moreover the homotopy $(H_t)_{t\in [0,1]}$ will have the property that $H_{(x,t)}:D\to A(x)$ is an isomorphism for all $x\in X$ and $t>0$. We prove this by induction on numbers $n$ with the property
that there are two closed covers of $X$, $W_1,..., W_n$ and $Y_1,...,Y_n$  such that $Y_i$ contained in the interior of $W_i$ and
 $A(W_i)\cong C(W_i)\otimes D$ for $1\leq i \leq n$.
First we observe that the case $n=1$ follows from Proposition~\ref{101}.
Let us now pass from $n-1$ to $n$. Given two covers as above, there is yet another closed cover $V_1,...,V_n$ of $X$ such that
$V_i$ is a neighborhood of $Y_i$ and $W_i$ is a neighborhood of $V_i$ for all $1\leq i \leq n$.
Set $Y=\cup_{i=1}^{n-1}Y_i$, $V=\cup_{i=1}^{n-1}V_i$ and $W=\cup_{i=1}^{n-1}W_i$.
By the inductive hypothesis applied to $A(V)$, and the covers $V_1,...,V_{n-1}$ and $ W_1\cap V,..., W_{n-1}\cap V$ there is a homotopy $h:D\to A(V)\otimes C[0,1]$ such that $h_{(x,0)}=\varphi_x$ 
and $h_{(x,t)}:D\to A(x)$ is an isomorphism for all $(x,t)\in V\times (0,1]$. Fix a trivialization 
$\nu:A(Y_{n+1})\to C(Y_{n+1})\otimes D$. Define a continuous map 
$\alpha:(V\cap Y_{n+1})\times [0,1] \cup Y_{n+1}\times \{0\}\to \End(D)$ by setting
$\alpha_{(x,t)}=\nu_x h_{(x,t)}$ if $(x,t)\in (V\cap Y_{n+1})\times [0,1]$ and
$\alpha_{(x,0)}=\nu_x \varphi_x$ if $x\in Y_{n+1}$.
Since $V\cap Y_{n+1}$ is a neighborhood of $Y\cap Y_{n+1}$ in $Y_{n+1}$ and since $\nu_x\varphi_x\in \End(D)^*$ for all $x\in Y_{n+1}$, by Lemma~\ref{102} there is a continuous map
$\eta:Y_{n+1}\times [0,1]\to \End(D)^*$ which extends the restriction of $\alpha$ to 
$(Y\cap Y_{n+1})\times [0,1] \cup Y_{n+1}\times \{0\}$.
We conclude the construction of the desired homotopy by defining $H:D\to A(X)\otimes C[0,1]$ by $H_{(x,t)}=h_{(x,t)}$ for
$(x,t)\in Y\times [0,1]$ and $H_{(x,t)}=\nu_x^{-1}\eta_{(x,t)}$ for
$(x,t)\in Y_{n+1}\times [0,1]$. 
\end{proof}
\begin{lemma}\label{local-triviality-CX-semiprojective}
  Let $D$ be a $KK$-semiprojective
stable Kirchberg algebra. Let  $X$ be a compact metrizable space and $Y,$ $Z$ be closed subsets of $X$ such that
$X=Y\cup Z$.
     Suppose that $\gamma:D \to C(Y\cap Z)\ot D$ is a
   full
    $*$-\ho\ which admits a lifting to
 a full
    $*$-\ho\ $\alpha:D \to C(Y)\ot D$
    such that $\alpha_x\in \End(D)^*$ for all $x\in Y$.
    Then the pullback $C(Y)\ot D\oplus_{\pi_{Y\cap Z},\widetilde{\gamma}\pi_{Y\cap Z}} C(Z)\ot D$ is isomorphic to
    $C(X)\ot D$.
\end{lemma}
\begin{proof} By
Prop.~\ref{101} there is a  $*$-homomorphism $\eta:D\to C(Y)\ot D$ such that
$\eta_x=\gamma_x$ for $x \in Y\cap Z$ and such that $\eta_x \in \Aut(D)$ for
$x\in Y\setminus Z$. Using the short five lemma one checks immediately that the triplet $(\widetilde\eta,\widetilde\gamma,
 \mathrm{id}_{C(Z)\ot D})$ defines a $C(X)$-linear
isomorphism: \newline $C(X)\ot D=C(Y)\ot D\oplus_{\pi_{Y\cap
Z},\pi_{Y\cap Z}} C(Z)\ot D\to C(Y)\ot D\oplus_{\pi_{Y\cap
Z},\widetilde{\gamma}\pi_{Y\cap Z}} C(Z)\ot D$.
\end{proof}
\begin{lemma}\label{local-triviality}
  Let $D$ be a $KK$-semiprojective stable Kirchberg
algebra. Let $Y$, $Z$ and $Z'$ be closed subsets of a compact metrizable space
$X$ such that $Z'$ is a neighborhood of $Z$ and $X=Y\cup Z$.
  Let $B$ be a  $C(Y)$-algebra locally
    isomorphic to $C(Y)\otimes D$
    and let $E$ be a  $C(Z')$-algebra locally
    isomorphic to $C(Z')\otimes D.$
     Let $\alpha: E(Y\cap Z')\to B(Y\cap Z')$  be
 a $*$-monomorphism of $C(Y\cap Z')$-algebras such that
$KK(\alpha_x)\in KK(E(x),B(x))^{-1}$ for all $x \in Y\cap Z'$. If
$\gamma=\alpha_{Y\cap Z}$,
   then $B(Y)\oplus_{\pi_{\SSS Y\cap Z},\gamma\,\pi_{\SSS Y\cap Z}}E(Z)$ is
locally isomorphic to $C(X)\ot D$.
\end{lemma}
\begin{proof} Since we are dealing with a local property,
we may assume that $B=C(Y)\ot D$ and $E=C(Z')\ot D$. To simplify notation we
let $\pi$ stand for both $\pi^Y_{Y\cap Z}$ and $\pi^Z_{Y\cap Z}$ in the sequel.
Let us denote by $H$ the $C(X)$-algebra $C(Y)\ot D
\oplus_{\pi,\gamma\pi}C(Z)\ot D.$
  We must show that $H$ is locally trivial.
 Let $x\in X$. If $x\notin  Z$, then
there is a closed neighborhood $V$ of $x$ which does not intersect $Z$, and
hence the restriction of $H$ to $V$ is isomorphic to $C(V)\ot D$, as it follows
immediately from the definition of $H$.
 It remains to consider the case
when $x\in Z$. Now $Z'$ is a closed  neighborhood of $x$ in $X$ and
 the
restriction of $H$ to $Z'$ is isomorphic to $C(Y\cap Z')\ot
D\oplus_{\pi,\gamma\pi}C(Z)\ot D$. Since $\gamma:Y\cap Z \to \End(D)^*$ admits
a continuous extension $\alpha:Y\cap Z' \to \End(D)^*$, it follows that $H(Z')$
is isomorphic to $C(Z')\ot D$ by
Lemma~\ref{local-triviality-CX-semiprojective}.
\end{proof}

\begin{proposition}\label{A(psi)-is-lt}
 Let $X$, $A$, $D$ and $\sigma$ be as in Proposition~\ref{approximation-fiber-is-D}.
 For any finite subset $\fset$ of $A$ and any $\ep>0$
 there is a  $C(X)$-algebra $B$ which is locally isomorphic to $C(X)\ot
D$ and there exists a  $C(X)$-linear
 $*$-monomorphism $\eta:B\to A$  such that $\fset\subset_\ep \eta(B)$ and
 $KK(\eta_x)\in KK(B(x),A(x))^{-1}$ for all $x\in X$.
\end{proposition}
\begin{proof} Let $\psi_k:E_k=C(Y_k)\ot D\to A(Y_k)$, $k=0,...,n$
 be as in the
conclusion of Proposition~\ref{approximation-fiber-is-D},  strengthen as in
Remark~\ref{polys}. Therefore we may assume that there is another $n$-fibered
$\{D\}$-monomorphism $(\psi'_0,...,\psi'_n)$ into $A$ such that
$\psi'_k:C(Y'_k)\ot D\to A(Y_k')$,  $Y_k'$ is a closed neighborhood of $Y_i$,
and $\pi_{Y_k}\psi'_k=\psi_k$, $k=0,...,n.$ Let $X_k$, $B_k$, $\eta_k$ and
$\gamma_k$ be as in Definition~\ref{stack}. $B_0$ and $\eta_0$ satisfy the
conclusion of the proposition, except that we need to prove that $B_0$ is
locally isomorphic to $C(X)\ot D$. We  prove by induction on decreasing $k$
that the $C(X_k)$-algebras $B_k$ are locally trivial. Indeed $B_n=C(X_n)\ot D$
and assuming that $B_k$ is locally trivial, it follows by
Lemma~\ref{local-triviality} that $B_{k-1}$ is locally trivial, since by
\eqref{eq-recurrence-B(k)}
\begin{equation*}\label{eq-recurrence-B(k)*}
B_{k-1}\cong B_k\oplus_{\pi\eta_k,\pi\psi_{k-1}} E_{k-1}\cong
B_k\oplus_{\pi,\gamma_k\pi} E_{k-1},\quad(\,\, \pi=\pi_{X_k\cap Y_{k-1}})
\end{equation*}
and $\gamma_k: E_{k-1}(X_k\cap Y_{k-1}) \to B_k(X_k\cap Y_{k-1})$,
$(\gamma_k)_x=(\eta_k)_x^{-1}(\psi_{k-1})_x$, extends to a $*$-\mo\ $\alpha:
E_{k-1}(X_k\cap Y'_{k-1}) \to B_k(X_k\cap Y'_{k-1})$,
$\alpha_x=(\eta_k)_x^{-1}(\psi'_{k-1})_x$
 and $KK(\alpha_x)$ is a KK-equivalence since both
$KK((\eta_k)_x)$ and $KK((\psi_{k-1})_x)$ are KK-equivalences.
\end{proof}

\section{When is a $C(X)$-algebra locally trivial}\label{LTC}

 In this section we prove  Theorems~\ref{cuntz-algebras-intro}
--~\ref{thm-lista-local-triviala}
 and   some of their consequences.
\vskip 4pt \emph{Proof of Theorem~\ref{stable-trivial-C(X)-algebras}.}

\begin{proof}Let $X$ denote the primitive spectrum of $A$. Then $A$ is a continuous
$C(X)$-algebra and its fibers are stable Kirchberg algebras (see
\cite[2.2.2]{BK}). Since $A$ is separable, $X$ is metrizable by
Lemma~\ref{X-is-2nd}.
  By Proposition~\ref{A(psi)-is-lt} there is a sequence of
 $C(X)$-algebras $(A_k)_{k=1}^\infty$ locally isomorphic to $C(X)\ot D$
  and a sequence of
 $C(X)$-linear $*$-monomorphisms
$(\eta_k:A_k\to A)_{k=1}^\infty$, such that $KK(\eta_k)_x$ is a KK-equivalence
for each $x \in X$ and $(\eta_k(A_k))_{k=1}^\infty$ is an exhaustive sequence
of $C(X)$-subalgebras of $A$. Since $D$ is weakly semiprojective and KK-stable,
after passing to a subsequence of $(A_k)$ if necessary, we find a sequence
$(\sigma_k)_{k=1}^\infty$, $\sigma_k\in KK(D,A_k)$ such that
$KK(\eta_k)\sigma_k=\sigma$ for all $k\geq 1$. Since both $KK(\eta_k)_x$ and
$\sigma_x$ are KK-equivalences, we deduce that $(\sigma_k)_x\in
KK(D,A_k(x))^{-1}$ for all $x\in X$. By
Proposition~\ref{loc-triv+section=triv}, for each $k\geq 1$ there is an
isomorphism of $C(X)$-algebras $\varphi_k:C(X)\ot D \to A_k$ such that
$KK(\varphi_k)=\sigma_k$. Therefore if we set $\theta_k=\eta_k\varphi_k$, then
$\theta_k$ is a  $C(X)$-linear $*$-monomorphism from $B\stackrel{def}=C(X)\ot
D$ to $A$ such that $KK(\theta_k)=\sigma$  and $(\theta_k(B))_{k=1}^\infty$ is
an exhaustive sequence of $C(X)$-subalgebras of $A$.
 Using again the weak semiprojectivity
and the KK-stability of $D$, and Lemma~\ref{fullness}, after passing to a
subsequence of $(\theta_k)_{k=1}^\infty$, we construct a sequence of finite
sets $\fset_k\subset B$ and a sequence of $C(X)$-linear  $*$-\mos\ $\mu_k:B\to
B$ such that
\begin{itemize}
      \item[(i)] $KK(\theta_{k+1}\mu_k)= KK(\theta_k)$ for all
      $k\geq 1$,
        \item[(ii)] $\|\theta_{k+1}\mu_k(a)-\theta_{k}(a)\|<2^{-k}$ for
        all $a \in \fset_k$ and all $k\geq 1$,
        \item[(iii)] $\mu_k(\fset_k)\subset \fset_{k+1}$ for all $k\geq 1$,
        \item[(iv)] $\bigcup_{j=k+1}^\infty\,
        \big(\mu_{j-1}\circ\cdots\circ\mu_{k}\big)^{-1}(\fset_j)$ is dense
        in $B$  and
        $\bigcup_{j=k}^\infty\,\theta_{j}(\fset_j)$ is
        dense in $A$ for all $k\geq 1$.
    \end{itemize}
    Arguing as in the proof of \cite[Prop.~2.3.2]{Ror:encyclopedia},
one verifies that $$\Delta_k(a)=\lim_{j \to
\infty}\theta_j\circ\big(\mu_{j-1}\circ\cdots\circ\mu_{k}\big) (a)$$ defines a
sequence of $*$-\mos\  $\Delta_k:B \to A$  such that
$\Delta_{k+1}\mu_k=\Delta_k$ and the induced map   $\Delta:\varinjlim_k
(B,\mu_k)\to A$ is an isomorphism of $C(X)$-algebras. Let us show that
$\varinjlim_k (B,\mu_k)$ is isomorphic to $B$. To this purpose, in view of
Elliott's intertwining argument, it suffices to show that each map $\mu_k$ is
approximately unitarily equivalent to a $C(X)$-linear automorphism of $B$.
Since $KK(\theta_k)=\sigma$, we deduce from (i) that $KK((\mu_k)_x)=KK(\id_D)$
for all $x\in X$. By Proposition~\ref{Pb}, this property implies that each map
$\mu_k$ is approximately unitarily equivalent to a $C(X)$-linear automorphism
of $B$. Therefore there is an isomorphism of $C(X)$-algebras $\Delta:B\to A$.
Let us show that we can arrange that $KK(\Delta|_D)=\sigma$. By
Theorem~\ref{Phillips-basic}, there is a full $*$-\ho\ $\alpha:D\to B$ such
that $KK(\alpha)=KK(\Delta^{-1})\sigma$.
 Since
$KK(\Delta_x^{-1})\sigma_x\in KK(D,D)^{-1}$, by Proposition~\ref{Pb} there is
$\Phi_1:D \to C(X)\ot D$ such that $\widetilde\Phi_1\in \Aut_{C(X)}(B)$ and
$KK(\Phi_1)=KK(\Delta^{-1})\sigma$. Then $\Phi=\Delta\widetilde\Phi_1:B\to A$
is an isomorphism such that $KK(\Phi|_D)=KK(\Delta\Phi_1)=\sigma$.
\end{proof}
Dixmier and Douady \cite{Dix:C*}  proved that a continuous field with
fibers $\mathcal{K}$ over a finite dimensional locally compact Hausdorff space
is locally trivial if and only it verifies Fell's condition, i.e. for each
$x_0\in X$ there is a continuous section $a$ of the field  such that $a(x)$ is
a rank one projection for each $x$ in a neighborhood of $x_0$. We have a
analogous result:
\begin{corollary}\label{cor-felie}
 Let $A$ be a  separable C*-algebra whose primitive
 spectrum $X$ is  Hausdorff and of finite dimension.
 Suppose that for each $x\in X$,
 $A(x)$  is $KK$-semiprojective,
 nuclear, purely infinite and stable.
Then $A$ is locally trivial if and only if for each $x\in X$
   there exist a closed
 neighborhood $V$ of $x$, a  Kirchberg algebra $D$ and
  $\sigma\in KK(D,A(V))$ such that
 $\sigma_v\in KK( D, A(v))^{-1}$ for each $v\in V$.
\end{corollary}
\begin{proof}One applies Theorem~\ref{stable-trivial-C(X)-algebras}
 for
$D\ot \mathcal{K}$ and $A(V)$.
\end{proof}
\begin{proposition}\label{prop:1_dim_cone}
 Let $\psi$ be a full endomorphism of a Kirchberg algebra $D$. 
If $D$ is unital we assume that $\psi(1)=1$ as well. Then the  continuous $C[0,1]$-algebra
$E=\{f\in  C[0,1]\ot D: f(0)\in\psi(D)\}$ is locally trivial if and only if $\psi$ is homotopic
to an automorphism of $D$.
\end{proposition}
\begin{proof} Suppose that $E$ is
 trivial on some neighborhood of $0$. Thus there is $s\in (0,1]$ and an isomorphism
 $\theta:C[0,s]\ot D\stackrel{\simeq}\to E[0,s]$. Since $E[0,s]\subset C[0,s]\ot D$,
there is a continuous path $(\theta_t)_{t\in [0,s]}$ in $\End(D)$ such that
 $\theta_t\in \Aut(D)$ for $0< t\leq s$ and $\theta_0(D)=\psi(D)$. Set
$\beta=\theta_0^{-1}\psi\in\Aut(D)$. Then $\psi$ is homotopic to an automorphism via the path $(\theta_t\beta)_{t\in [0,s]}$. Conversely, if $\psi$ is homotopic to an automorphism $\alpha$,
then by Theorem~\ref{Phillips-basic}  there is
a continuous path $(u_t)_{t\in (0,1]}$ of unitaries in $D^+$ such that $\lim_{t\to 0}\|\psi(d)-u_t \alpha(d)u_t^*\|=0$
for all $d\in D$. The path $(\theta_t)_{t\in [0,1]}$ defined by $\theta_0=\psi$ and $\theta_t=u_t \alpha u_t^*$ for $t\in (0,1]$ induces a $C[0,1]$-linear $*$-endomorphism
of $C[0,1]\ot D$ which maps injectively $C[0,1]\ot D$ onto $E$.

\end{proof}

\vskip 4pt \emph{Proof of Theorem~\ref{o2k}.} \begin{proof}For the first part
we apply Theorem~\ref{stable-trivial-C(X)-algebras} for
$D=\OOO_2\ot\mathcal{K}$ and $\sigma=0$. For the second part we assert that if
$D$ is a  Kirchberg such that all continuous $C[0,1]$-algebras with fibers
isomorphic to $D$ are locally trivial then $D$ is stable and $KK(D,D)=0$. Thus
 $D$ is KK-equivalent to $\OOO_2$ and hence that $D\cong
\OOO_2\ot\mathcal{K}$ by \cite[Thm.~8.4.1]{Ror:encyclopedia}. The Kirchberg
algebra D is either unital or stable \cite[Prop.~4.1.3]{Ror:encyclopedia}. Let
$\psi:D\to D$ be a $*$-\mo\  such that $KK(\psi)=0$ and such that
$\psi(1_D)<1_D$ if $D$ is unital. By Proposition~\ref{prop:1_dim_cone} $\psi$ is homotopic to an automorphism of $\theta$ of $D$.
 Therefore $D$ must be nonunital (and
hence stable), since otherwise $1_D$ would be homotopic to its proper
subprojection $\psi(1_D)$. Moreover $KK(\theta)=KK(\psi)=0$ and hence
$KK(D,D)=0$ since $\theta$ is an automorphism.
\end{proof}

 We turn now to unital $C(X)$-algebras.
\begin{theorem}\label{trivial-C(X)-algebras}  Let $A$ be a separable unital
$C(X)$-algebra over a finite dimensional compact Hausdorff space $X$.
  Suppose that  each fiber
 $A(x)$  is
  nuclear simple and purely infinite.
   Then $A$ is isomorphic to $C(X)\ot D$,
   for some   $KK$-semiprojective unital Kirchberg algebra $D$,
if and only if there is $\sigma \in KK(D,A)$ such that
 $K_0(\sigma)[1_D]=[1_A]$
and $\sigma_x \in KK(D, A(x))^{-1}$ for all $x \in X$. For any such $\sigma$
there is an  isomorphism of $C(X)$-algebras $\Phi:C(X)\ot D\to A$ such that
$KK(\Phi|_D)=\sigma$.
 \end{theorem}
\begin{proof} We verify the nontrivial implication. $X$ is metrizable by
Lemma~\ref{X-is-2nd}. $A$ is a continuous $C(X)$-algebra by
Lemma~\ref{lemma-auto-continuity}.
 By
Theorem~\ref{stable-trivial-C(X)-algebras}, there is an isomorphism
$\Phi:C(X)\ot D\ot \mathcal{K}\to A\ot \mathcal{K}$ such that
$KK(\Phi)=\sigma$. Since $K_0(\sigma)[1_D]=[1_A]$, and since $A\ot \mathcal{K}$
contains a full properly infinite projection, we may arrange that
$\Phi(1_{C(X)\ot D}\ot e_{11})=1_A\ot e_{11}$ after conjugating $\Phi$ by some
unitary $u\in M(A\ot \mathcal{K})$. Then $\varphi=\Phi|_{C(X)\ot D\ot e_{11}}$
satisfies the conclusion of the theorem.
\end{proof}
 \emph{Proof of Theorem~\ref{A}.}

\begin{proof}Let $D$ be a KK-semiprojective unital Kirchberg algebra $D$
such that every unital $*$-endomorphism of $D$ is a KK-equivalence. Suppose
that $A$ is a separable unital $C(X)$-algebra over a finite dimensional
 compact Hausdorff space the fibers of which are isomorphic to $D$.
 We shall prove that $A$ is locally trivial. By Theorem~\ref{trivial-C(X)-algebras},
  it suffices
 to show that each point $x_0\in X$ has a closed neighborhood $V$ for which there is 
$\sigma\in KK(D,A(V))$ such that  $K_0(\sigma)[1_D]=[1_{A(V)}]$ and $\sigma_x\in KK(D,A(x))^{-1}$ for all $x\in V$.

 Let $(V_n)_{n=1}^\infty$ be a decreasing sequence of closed neighborhoods of $x_0$
whose intersection is $\{x_0\}$.
 Then $A(x_0)\cong\varinjlim \,A(V_n).$  By assumption, there is an
 isomorphism $\eta:D\to A(x_0)$.
 Since $D$ is KK-semiprojective, there is $m\geq 1$ such that $KK(\eta)$ lifts to some
 $\sigma\in KK(D,A(V_m))$ such that $K_0(\sigma)[1_D]=[1_{A(V_m)}]$.
Let $x\in V_m$.
By assumption, there is an isomorphism $\varphi:A(x)\rightarrow D$.
 The $K_0$-morphism induced by $KK(\varphi)\sigma_x$ maps $[1_D]$ to itself.
 By Theorem~\ref{Phillips-basic} there is a unital $*$-\ho\
 $\psi:D\to D$ such that
 $KK(\psi)=KK(\varphi)\sigma_x$. By assumption we must have $KK(\psi)\in KK(D,D)^{-1}$
  and hence $\sigma_x\in KK(D,A(x))^{-1}$
 since $\varphi$ is an isomorphism.
 Therefore $A(V_m)\cong C(V_m)\ot D$ by Theorem~\ref{trivial-C(X)-algebras}.

Conversely, let us assume that
 all separable unital continuous $C[0,1]$-algebras
 with fibers isomorphic to $D$
 are locally trivial. 
Let $\psi$ be any unital $*$-endomorphism of $D$. 
By Proposition~\ref{prop:1_dim_cone} $\psi$ is homotopic to an automorphism of $D$ and 
    hence $KK(\psi)$ is invertible.
\end{proof}

\vskip 4pt \emph{Proof of Theorem~\ref{cuntz-algebras-intro}}
\begin{proof}
 Let $A$ be as in
Theorem~\ref{cuntz-algebras-intro} and let $n\in \{2,3,...\}\cup\{\infty\}$. It
is known that $\OOO_n$ satisfies the UCT. Moreover $K_0(\OOO_n)$ is generated
by $[1_{\OOO_n}]$ and $K_1(\OOO_n)=0$. Therefore any unital $*$-endomorphism of
$\OOO_n$ is a KK-equivalence. It follows that $A$ is locally trivial by
Theorem~\ref{A}. Suppose now that $n=2$. Since
$KK(\OOO_2,\OOO_2)=KK(\OOO_2,A)=0$, we may apply Theorem~\ref{A} with
$\sigma=0$ and obtain that $A\cong C(X)\ot \OOO_2$. Suppose now that
$n=\infty$. Let us define $\theta:K_0(\OOO_\infty)\to K_0(A)$ by
$\theta(k[1_{\OOO_\infty}])=k[1_A]$, $k\in \ZZZ$. Since $\OOO_\infty$ satisfies
the UCT, $\theta$ lifts to some element $\sigma\in KK(\OOO_\infty,A)$.  By
Theorem~\ref{A} it follows that $A\cong C(X)\ot \OOO_\infty$. Finally let us
consider the case $n\in \{3,4,...\}$. Then $K_0(\OOO_n)=\mathbb{Z}/(n-1)$.
Since $\OOO_n$ satisfies the UCT, the existence of an element $\sigma\in
KK(\OOO_n,A)$ such that  $K_0(\sigma)[1_{\OOO_n}]=[1_A]$ is equivalent to the
existence of a
 morphism of groups
 $\theta:\mathbb{Z}/(n-1) \to K_0(A)$ such that $\theta(\bar{1})=[1_A]$.
 This is equivalent to requiring  that $(n-1)[1_A]=0$.
\end{proof}

 As a corollary of Theorem~\ref{cuntz-algebras-intro} we have that
$[X,\Aut(\OOO_\infty)]$ reduces to a point.  The homotopy groups of the
endomorphisms of the stable Cuntz-Krieger algebras were computed in
\cite{Cuntz:homotopy-Aut}. Let $v_1,\dots,v_n$ be the canonical generators of
$\OOO_n$, $2 \leq n <\infty$.
\begin{theorem}\label{nontrivialO(n)bundle}
For any   compact metrizable space $X$ there is a bijection
$[X,\Aut(\OOO_n)]\to K_1(C(X)\ot \OOO_n)$.
    The $k^{th}$-homotopy group $\pi_{k}(\mathrm{Aut}(\OOO_n))$ is isomorphic to
    $\ZZZ/(n-1)$
    if $k $ is odd and it vanishes if $k$ is even.
    In particular
    $\pi_{1}(\mathrm{Aut}(\OOO_n))$ is generated by the class
    of the canonical action of $\mathbb{T}$
     on $\OOO_n$,
    $\lambda_z(v_i)=zv_i$.
\end{theorem}
\begin{proof} Since $\OOO_n$
satisfies the UCT, we deduce that $\End(\OOO_n)^*=\End(\OOO_n)$. An immediate
application of Proposition~\ref{Pb} shows that the natural map
$\Aut(\OOO_n)\hookrightarrow \End(\OOO_n)$ induces an isomorphism of groups
$[X,\Aut(\OOO_n)]\cong[X,\End(\OOO_n)]$.  Let $\iota:\OOO_n \to C(X)\ot \OOO_n$
be defined by $\iota(v_i)=1_{C(X)}\ot v_i$, $i=1,...,n$. The map $\psi\mapsto
u(\psi)=\psi(v_1)\iota(v_1)^*+\cdots+\psi(v_n)\iota(v_n)^*$ is known to be a
homeomorphism from $\mathrm{Hom}(\OOO_n, C(X)\ot\OOO_n)$ to the unitary group
of $C(X)\ot \OOO_n$. Its inverse maps a unitary $w$ to the $*$-\ho\ $\psi$
uniquely defined by $\psi(v_i)=w\iota(v_i)$, $i=1,...,n$. Therefore
\[[X,\Aut(\OOO_n)]\cong [X,\End(\OOO_n)]\cong  \pi_0(U(C(X)\ot \OOO_n))\cong
K_1(C(X)\ot \OOO_n).\] The last isomorphism holds since
 $\pi_0(U(B))\cong K_1(B)$ if $B\cong B\ot \OOO_\infty$, by
\cite[Lemma~2.1.7]{Phi:class}.  One verifies easily that if
$\varphi\in\mathrm{Hom}(\OOO_n, C(X)\ot\OOO_n)$, then
$u(\widetilde\psi\varphi)=\widetilde\psi(u(\varphi))u(\psi)$. Therefore the
bijection $\chi:[X,\End(D)] \to K_1(C(X)\ot \OOO_n)$ is an isomorphism of
groups whenever $K_1(\widetilde\psi)=\id$ for all $\psi\in \mathrm{Hom}(\OOO_n,
C(X)\ot\OOO_n)$. Using the $C(X)$-linearity of $\widetilde\psi$ one observes
that this holds  if the $n-1$ torsion of $ K_0(C(X))$ reduces to $\{0\}$, since
in that case the map $K_1(C(X))\to K_1(C(X)\ot\OOO_n)$ is surjective by the
K\"{u}nneth formula.
\end{proof}
\begin{corollary}\label{2-sphere} Let $X$ be a finite dimensional compact
metrizable space. The isomorphism classes of unital separable $C(SX)$-algebras
    with all fibers isomorphic to $\OOO_n$  are
     parameterized by $K_1(C(X)\ot \OOO_n)$.
\end{corollary}
\begin{proof} This follows from Theorems~\ref{cuntz-algebras-intro}
 and ~\ref{nontrivialO(n)bundle}, since the locally trivial
 principal $H$-bundles over  $SX=X\times[0,1]/X\times \{0,1\}$
  are parameterized by the homotopy classes $[X,H]$ if $H$
 is a path connected group \cite[Cor.~8.4]{Hus:fibre}. Here we take
 $H=\Aut(\OOO_n)$.
\end{proof}
Examples of nontrivial unital $C(X)$-algebras with fiber $\OOO_n$ over a
$2m$-sphere arising from vector bundles were exhibited in 
\cite{Vasselli:polish-bundles}, see also \cite{Vasselli:jfa-bundles}.

We need some preparation for the proof of
Theorem~\ref{thm-lista-local-triviala}.
 Let $G$ be a group, let $g\in G$ and  set
$\End(G,g)=\{\alpha\in \End(G):\alpha(g)=g\}$.
 The pair $(G,g)$ is called
\emph{\weakr} if $\End(G,g)\subset \Aut(G)$ and \emph{rigid} if $\End(G,g)=\{ \id_G\}$.

\begin{theorem}\label{all-wr}
If $G$ is a finitely generated abelian group, then $(G,g)$ is
 \emph{\weakr} if and only if $(G,g)$ is isomorphic to one of the pointed groups from the list
$\mathcal{G}$ of Theorem~\ref{thm-lista-local-triviala}.
\end{theorem}
\begin{proof} First we make a number of remarks.

 (1)   $(G,g)$ is \weakr\ if and only if $(G,\alpha(g))$ is \weakr\ for some (or any) $\alpha\in \mathrm{Aut}(G)$. Indeed if $\beta \in \End(G,g)$ then $\alpha\beta\alpha^{-1}\in \End(G,\alpha(g))$.

(2) By considering the zero endomorphism of $G$ we see that if $(G,g)$ is \weakr\ and $G\neq 0$ then $g\neq 0$. 

(3) If $(G\oplus H,g\oplus h)$ is \weakr, then so are
$(G,g)$ and $(H,h)$. 

(4) Let us observe that $(\mathbb{Z}^2,g)$ is not \weakr\ for any $g$. Indeed,
if $g=(a,b)\neq 0$, then the matrix $\begin{pmatrix} 1+b^2 & -ab\\ -ab &
1+a^2\end{pmatrix}$ defines an endomorphism $\alpha$ of $\mathbb{Z}^2$ such
that $\alpha(g)=g$, but $\alpha$ is not invertible since
$\mathrm{det}(\alpha)=1+a^2+b^2>1$.

(5)  Let $p$ be a prime and let $1\leq e_1 \leq e_2$, $0\leq s_1 < e_1$, $0\leq s_2 < e_2$ be  integers. If $(G,g)=(\Z/p^{e_1}\oplus\Z/p^{e_2},p^{s_1}\oplus p^{s_2})$ is \weakr\ then $0<s_2-s_1<e_2-e_1$. Indeed if $s_1\geq s_2$ then the matrix 
$\begin{pmatrix} 0 & p^{s_1-s_2}\\ 0 &
1\end{pmatrix}$ induces a noninjective endomorphism of $(G,g)$. Also if $s_1<s_2$ and
 $s_2-s_1\geq e_2-e_1$
then $p^{e_1}\bar{b}=0$ in  $\Z/p^{e_2}$, where $b=p^{s_2-s_1}$ and so the matrix $\begin{pmatrix} 1 & 0\\b  &
0\end{pmatrix}$ induces a well-defined noninjective endomorphism of $(G,g)$.

(6) Let $p$ be a prime and let $1\leq k$, $0\leq s < e$ be integers.
Suppose that $(\Z\oplus \Z/p^e,k\oplus p^s)$ is \weakr. Then $k$ is divisible by $p^{s+1}$.
Indeed, seeking a contradiction suppose that $k$ can be written as
 $k=p^t c$ where $0\leq t \leq s$ and $c$ are integers such that $c$ is not divisible by $p$.
Let $d$ be an integer such that $dc-1$ is divisible by $p^{e}$. Then the matrix
$\begin{pmatrix} 1 & 0\\dp^{s-t}  &
0\end{pmatrix}$ induces a  noninjective endomorphism of $(\Z\oplus \Z/p^{e},k\oplus p^{s})$.

Suppose now that $(G,g)$ is \weakr. We shall show that $(G,g)$ is isomorphic to one of the pointed groups from the list $\mathcal{G}$. Since $G$ is abelian and
 finitely generated it decomposes  as a direct sum
of its primary components
\begin{equation}\label{eqn:ab_structure}
 G\cong \Z^r \oplus G(p_1)\oplus \cdots \oplus G(p_m)
\end{equation}
where $p_i$ are distinct prime numbers. 
Each primary component $G(p_i)$ is of the form
\begin{equation}\label{eqn:prime_structure}
G(p_i)=\Z/p_i^{e_{i\, 1}}\oplus \cdots\oplus Z/p_i^{e_{i \,n(i)}}
\end{equation}
where $1\leq e_{i\,1}\leq \cdots \leq e_{i\, n(i)}$ are positive integers.
Corresponding to the decomposition \eqref{eqn:ab_structure} we write the base point $g=g_0\oplus g_1\oplus ...\oplus g_m$
with $g_0\in \Z^r$ and $g_i\in G(p_i)$ for $i\geq 1$.
If $g_{ij}$ is the component of $g_i$ in $\Z/p^{e_{i j}}$, then it follows from 
(1), (2) and (3) that we may assume that $g_{i j}=p^{s_{i j}}$ for some integer 
$0\leq s_{ij}<e_{ij}$.
Using (3) and (4) we deduce that $r=1$ in \eqref{eqn:ab_structure} and that
$g_0=k \neq 0$ by (2). We may assume that $k\geq 1$ by (1).
Then using (3) and (5) we deduce that for each $1\leq i \leq m$, 
$0<s_{i\, j+1}-s_{i\, j}<e_{i\, j+1}-e_{i\, j}$ for $1\leq j < n(i)$.
Finally, from (3) and (6) we see that $k$ is divisible by the product 
$p_1^{s_{1\, n(1)}}\cdots p_m^{s_{m\, n(m)}}$. Therefore $(G,g)$ is isomorphic 
to one of the pointed groups on the list $\mathcal{G}$.

Conversely, we shall prove that if $(G,g)$ belongs to the list $\mathcal{G}$ then $(G,g)$
is \weakr. This is obvious if $G$ is torsion free i.e for $(\{0\},0)$ 
and $(\Z,k)$ with $k\geq 1$. 

Let us consider the case when $G$ is a torsion group.
Since
\[\End\big(G(p_1)\oplus \cdots \oplus G(p_m), g_1\oplus \cdots \oplus g_m \big)
\cong \bigoplus_{i=1}^m\End\big(G(p_i),g_i)\]
it suffices to assume that $G$ is a $p$-group,
\[(G,g)=\big(\Z/p^{e_1}\oplus \cdots\oplus Z/p^{e_n}, p^{s_1}\oplus\cdots\oplus p^{s_n}\big)\]
with $0\leq s_i < e_i$ for $i=1,...,n$ and $0<s_{i+1}-s_i<e_{i+1}-e_i$ for  $1\leq i<n$.
 For each $0\leq i,j\leq n$ set
$e_{ij}=\max\{e_i-e_j,0\}$. It follows immediately 
that $s_i < e_{ij}+s_j$ for all $i\neq j$.
Let $\alpha\in \End(G,g)$. It is well-known  that $\alpha$ is induced by
a square matrix $A=[a_{ij}]\in M_n(\Z)$ with the property that each entry $a_{ij}$
is divisible by $p^{e_{ij}}$ and so $a_{ij}=p^{e_{ij}}b_{ij}$ for some $b_{ij}\in Z$, see \cite{Hillar:monthly}. 
Since $\alpha(g)=g$, we  have $\sum_{j=1}^n \bar{b}_{ij}p^{e_{ij}+s_j}=p^{s_i}$ in $\Z/p^{e_i}$
for all $0\leq i\leq n$. Since $e_{ij}+s_j>s_i$ for $i\neq j$ and $e_i>s_i$ we see that $b_{ii}-1$ must be divisible by $p$ for all $1\leq i \leq n$. Since $\mathrm{det}(A)$ is congruent to $b_{11}\cdots b_{nn}$ modulo $p$
it follows that  $\mathrm{det}(A)$ is not divisible by $p$ and so $\alpha\in \mathrm{Aut}(G)$  by \cite{Hillar:monthly}.

Finally consider the case when $(G,g)=\big(\Z\oplus G(p_1)\oplus\cdots \oplus G(p_m), k\oplus g_1\oplus\cdots \oplus g_m\big).$
If $\gamma \in \End(G,g)$ then there exist $\alpha_i \in \End (G(p_i),g_i)$ 
and $d_i\in G(p_i)$, $1\leq i \leq n$,
such that $\gamma (x_0\oplus x_1\oplus...\oplus x_m)=x_0\oplus (\alpha_1(x_1)+x_0 d_1)\oplus...\oplus(\alpha_m(x_m)+x_0 d_m)$.
Note that if each $\alpha_i$ is an automorphism then so is $\gamma$. Indeed, its inverse is
$\gamma^{-1}(x_0\oplus x_1\oplus...\oplus x_m)=x_0\oplus (\alpha_1^{-1}(x_1)+x_0 c_1)\oplus...\oplus(\alpha_m(x_m)^{-1}+x_0 c_m)$,
where $c_i=-\alpha_i^{-1}(d_i)$. Therefore it suffices to consider the case $m=1$, i.e.
\[(G,g)=\big(\Z\oplus\Z/p^{e_1}\oplus \cdots\oplus Z/p^{e_n},k\oplus p^{s_1}\oplus\cdots\oplus p^{s_n}\big),\]
and $(G,g)$ is on the list $\mathcal{G}$ (e). In particular $k=p^{s_n+1}\ell$ for some $\ell \in \mathbb{Z}$.
Let $\gamma \in \End(G,g)$. Then there exists $\alpha \in \End (G(p))$ and $d\in G(p)$
such that $\gamma(x_0\oplus x)=x_0 \oplus (\alpha(x)+x_0 d)$.
Just as above, $\alpha$ is induced by  a square matrix $A\in M_{n}(\Z)$ of the form 
 $A=[b_{ij}p^{e_{ij}}]\in M_{n}(\Z)$ with $b_{ij}\in \Z$, $e_{ij}=\max\{e_i-e_j,0\}$.
Since $\gamma(g)= g$ we have that 
$p^{s_n+1}\ell d_i+\sum_{j=1}^n \bar{b}_{ij}p^{e_{ij}+s_j}=p^{s_i}$ in $\Z/p^{e_i}$
for all $0\leq i\leq n$,  where the $d_i$ are the components of $d$. By reasoning as in the case when $G$ was a torsion group considered above, since  $s_{n}+1>s_i$ for all $1 \leq i \leq n$, $ e_{ij}+s_j>s_i$ for all $i\neq j$ and $e_i>s_i$, it follows again that each $b_{ii}-1$ is divisible by $p$ and that the endomorphism $\alpha$ of $G(p)$ induced by the matrix $A$ is an
automorphism. We conclude that $\gamma$ is an automorphism.
\end{proof}

\emph{Proof of Theorem~\ref{thm-lista-local-triviala}} \begin{proof} (ii) and (iii) Let $D$ be
a unital Kirchberg algebra such that $D$ satisfies the UCT and $K_*(D)$ is
finitely generated. Then $D$ is KK-semiprojective by Proposition~\ref{fg-uct}
and $KK(D,D)^{-1}=\{\alpha\in KK(D,D): K_*(\alpha)\,\,\text{is bijective}\}$.
In conjunction with Theorem~\ref{Phillips-basic}, this shows that all unital $*$-endomorphisms of $D$ are
KK-equivalences if and only if  both  $(K_0(D), [1_D])$ and $(K_1(D), 0)$
are \weakr. Equivalently, $K_1(D)=0$ and  $(K_0(D),[1_D])$ is \weakr. By
Theorem~\ref{all-wr}  $(K_0(D),[1_D])$ is \weakr\ if and only if
it isomorphic to one pointed groups from the list $\mathcal{G}$ of Theorem~\ref{thm-lista-local-triviala}.
 We conclude  the proof of (ii) and (iii) by applying Theorem~\ref{A}.

(i) By Theorem~\ref{cuntz-algebras-intro} both $\OOO_2$ and $\OOO_\infty$ 
have the automatic triviality property.
Conversely,  suppose that $D$ has the automatic triviality property,  where $D$ is a unital
Kirchberg algebra satisfying the UCT and such that $K_*(D)$  is finitely
generated.
   We shall prove that
   $D$ is isomorphic
 to either $\OOO_2$ or $\OOO_\infty$.

Let $Y$ be a finite connected CW-complex and let $\iota:D \to C(Y)\otimes D$ be the map $\iota(d)=1\otimes d$.
Let $[D,C(Y)\otimes D]$ denote the homotopy classes of unital $*$-homomorphisms from $D$ to $C(Y)\otimes D$.
By Theorem~\ref{Phillips-basic} the image of the map $\Delta:[D,C(Y)\otimes D]\to KK(D,C(Y)\otimes D)$ defined by
$[\varphi]\mapsto KK(\varphi)-KK(\iota)$ coincides with the kernel of the restriction morphism
$\rho: KK(D,C(Y)\otimes D)\to KK(\mathbb{C}1_D,C(Y)\otimes D)$.

 We claim that $\ker \rho$ must vanish
for all $Y$. Let $h \in\ker \rho$.
Then there is a unital $*$-homomorphism
$\varphi:D\to C(Y)\otimes D$ such that $\Delta[\varphi]=h$.
By Theorem~\ref{A}, each unital endomorphism of $D$ induces a KK-equivalence.
Therefore, by Proposition~\ref{Pb} there is a
$*$-\ho\ $\Phi:D\to C(Y)\otimes D$ such that
$\Phi_y\in \Aut(D)$ for all $y\in Y$ and $KK(\Phi)=KK(\varphi)$. Therefore $\Delta[\Phi]=KK(\Phi)-KK(\iota)=h$.
By hypothesis, the $\Aut(D)$-principal bundle constructed over the
suspension of $Y$ with characteristic map $y\mapsto \Phi_y$ is  trivial.
It follows then from  \cite[Thm.~8.2 p85]{Hus:fibre} that this map is homotopic
to the to the constant map $Y\to\Aut(D)$ which shrinks $Y$ to $\id_D$. This implies
that $\Phi$ is homotopic to $\iota$ and hence $h=0$.

Let us now observe that $\ker \rho$ contains subgroups isomorphic to $\mathrm{Hom}(K_1(D),K_1(D))$  and $\mathrm{Ext}(K_0(D),K_0(D))$ if $Y=\mathbb{T}$, since $D$ satisfies the UCT.
 It follows that both these groups  must vanish and so $K_1(D)=0$ and $K_0(D)$ is torsion free. 
On the other hand, 
 $(K_0(D),[1_D])$ is weakly rigid by the first part of the proof. Since $K_0(D)$ is torsion free we deduce from Theorem~\ref{all-wr}
that either
  $K_0(D)=0$ in which case $D\cong \OOO_2$ or that
  $(K_0(D),[1_D])\cong (\mathbb{Z},k)$, $k\geq 1$,
 in which case $D\cong
 M_{k}(\OOO_\infty)$
 by the classification
theorem of Kirchberg and Phillips.

 To conclude the proof,  it suffices to show
that $\ker \rho\neq 0$ if $D=M_k(\OOO_\infty)$, $k\geq 2$ and $Y$
is the two-dimensional space obtained by attaching a disk to a
circle by a degree-$k$ map. 
Since $K_0(C(Y)\ot \OOO_\infty)\cong
\mathbb{Z}\oplus \mathbb{Z}/k$ we can identify the map $\rho$ with the map
$\mathbb{Z}\oplus \mathbb{Z}/k \to \mathbb{Z}\oplus \mathbb{Z}/k$, $x\mapsto kx$
and so $\ker \rho\cong \mathbb{Z}/k\neq 0$ if $k\geq 2$.
\end{proof}

\emph{Added in proof.}
Some of the results from this paper are further developed in \cite{Dad:KKX-equiv}.
Theorem~\ref{stable-trivial-C(X)-algebras} was shown  to hold for all stable Kirchberg 
algebras $D$. 
The assumption that $X$ is finite dimensional is essential Theorem~\ref{cuntz-algebras-intro}.
Theorem~\ref{thm-lista-local-triviala} (ii) extends as follows: $\OOO_2$, $\OOO_\infty$ and $B\otimes \OOO_\infty$,
where $B$ is a unital UHF algebra of infinite type, are the only unital Kirchberg algebras
which satisfy the UCT and have the automatic triviality property.

\end{document}